\documentclass[11pt]{article}  % include bezier curves
 \pdfoutput=1
\usepackage{color}              % Need the color package
\usepackage{epsfig}
\usepackage{graphicx}
\usepackage{amssymb}
\usepackage{epsfig}
\usepackage{bm}

\newtheorem{Lemma}{Lemma}
\newtheorem{Proposition}[Lemma]{Proposition}

\newtheorem{Theorem}[Lemma]{Theorem}

\newtheorem{Corollary}[Lemma]{Corollary}

\newcommand{\ed}{\ \stackrel{d}{=} \ }
\newcommand{\cd}{\ \stackrel{d}{\rightarrow} \ }
\newcommand{\cp}{\ \stackrel{p}{\rightarrow} \ }
\newcommand{\FF}{\mathcal{F}}
\newcommand{\GG}{\mathcal{G}}
\newcommand{\NN}{\mathcal{N}}
\newcommand{\RR}{\mathcal{R}}

\newcommand{\WW}{\mathcal{W}}
\newcommand{\eps}{\varepsilon}
\newcommand{\bt}{\mathbf{t}}

\newcommand{\bT}{\mathbf{T}}

\newcommand{\var}{{\rm var}\ }

\newcommand{\sfrac}[2]{{\textstyle\frac{#1}{#2}}}
\newcommand{\qed}{\ \ \rule{1ex}{1ex}}

\newcommand{\len}{{\rm len}}

\newcommand{\Rbold}{{\mathbb{R}}}

\newcommand{\bpi}{\mbox{\boldmath$\pi$}}

\newcommand{\wL}{\widetilde{L}}
\newcommand{\wT}{\widetilde{T}}
\newcommand{\wW}{\widetilde{W}}

\newcommand{\proof}{{\bf Proof.\ }}
\renewcommand{\Pr}{{\mathbb{P}}}
\newcommand{\Ex}{{\mathbb{E}}}
\newcommand{\Reals}{{\mathbb{R}}}

\newcommand{\expec}{\mathbb{E}}
\newcommand{\prob}{\mathbb{P}}
\newcommand{\conde}{\mathbb{E}_\tau}
\newcommand{\condp}{\mathbb{P}_\tau}

\def\ind{{\rm 1\hspace{-0.90ex}1}}
\newcommand{\Lip}{{\mbox{\tiny{Lip}}}}
\begin{document}
\author{David J. Aldous\thanks{Research supported by N.S.F. Grant
DMS0704159}
\ and Shankar Bhamidi
\\
\\
       University of California\\
       Department of Statistics\\
        367 Evans Hall \# 3860\\
       Berkeley CA 94720-3860}

\title{Edge Flows in the Complete Random-Lengths Network}

\maketitle

\begin{abstract}
Consider the complete $n$-vertex graph whose edge-lengths are independent
exponentially distributed random variables.  Simultaneously for each pair of
vertices, put a constant flow between them along the shortest path.
Each edge gets some random total flow.
In the $n \to \infty$ limit we find explicitly the empirical distribution
of these edge-flows, suitably normalized.
\end{abstract}

\vspace{0.6in}


{\bf MSC2000 subject classification.}
60C05, 05C80, 90B15.

\newpage
\section{Introduction}
\label{sec-int}

Write {\em network} for an undirected graph whose edges $e$ have positive real
edge-lengths $\ell(e)$.
In a $n$-vertex connected network, the 
{\em distance} $D(i,j)$ between vertices $i$ and $j$ is the length of the shortest
route between them.  
Assuming generic edge-lengths, the shortest route is unique.
For each ordered (source, destination) pair of vertices
$(i,j)$, send flow of volume $1/n$ along the shortest route
from $i$ to $j$.  
(The normalization $1/n$ is arbitrary but convenient for (\ref{feD}) below).  
For each directed edge $e$ (i.e. an edge $e$ and a specified direction across $e$)
of the network, let $f(e)$ be the total
flow across the edge in that direction.  Note
\begin{equation}
n^{-1} \sum_{\mbox{\tiny{directed }} e}  f(e) \ell(e) 
 = n^{-2} \sum_i \sum_j D(i,j)
:= \bar{D} \label{feD}
\end{equation} 
where $\bar{D}$ is the average vertex-vertex distance.  

One can formulate a project to study the distribution of such edge-flows $f(e)$ in
different models of random $n$-vertex networks.  Such models include both
deterministic graphs to which random edge-lengths are assigned, and random graphs of
both the classical 
Erd\H{o}s - R\'{e}nyi  or random regular type \cite{bol85}
and the more recent {\em complex networks} types
\cite{AB02,durr-RGD,MD03,newman-survey}
again with real edge-lengths attached.
As (\ref{feD}) implies, this project is a refinement of the project of studying 
$\Ex \bar{D}_n$, so we envisage a model sufficiently tractable that we know
\begin{equation}
 \Ex \bar{D}_n = (1 + o(1)) \bar{d}_n  \label{Dn1}
\end{equation}
for some explicit $(\bar{d}_n)$.  

To set up some notation, return to the setting of a deterministic network.
Because we are using shortest-path routing, we expect edges-flows to be correlated with edge-lengths,
so let us study jointly edge-flows and edge-lengths by
considering the empirical measure $\psi^0$ which puts weight $1/n$ on each point
$(\ell(e),f(e))$:
\[ \psi^0(\cdot, \cdot) := \sfrac{1}{n} \sum_{\mbox{\tiny{directed }} e}  \ind
\{(\ell(e),f(e)) \in (\cdot, \cdot)\} . \]
So (\ref{feD}) becomes
\[ \int \int \ell u \ \psi^0(d\ell, du) = \bar{D} . \] 
So when short edge-lengths are order $1$ we should normalize edge-flows by
$\bar{d}_n$, that is consider the measure 
\begin{equation}
 \psi_n(\cdot, \cdot) := \sfrac{1}{n} \sum_{\mbox{\tiny{directed }} e}  \ind
\{(\ell(e),f(e)/\bar{d}_n) \in (\cdot, \cdot)\} .\label{psi-n} 
\end{equation}
For a random network, $\bar{D}_n$ is a random variable and $\psi_n(\cdot, \cdot)$ is
a random measure, related by  
\begin{equation}
 \int \int \ell y \ \psi_n(d\ell, dy) = \frac{\bar{D}_n}{\bar{d}_n} . 
 \label{ell-y}
\end{equation}
This notation is designed to suggest possible $n \to \infty$ limit behavior; that
the random 
measures $\psi_n$ converge to a non-random measure $\psi$ which by (\ref{Dn1}) and
under appropriate uniform integrability conditions must satisfy
\begin{equation}
 \int \int \ell y \ \psi(d\ell, dy) = 1. \label{psi-1}
\end{equation}
The purpose of this paper is to prove this result and identify $\psi$ in one
particular model, described in the next section.  
There is a fairly simple heuristic argument to identify $\psi$, shown in section \ref{sec-heur}. The
heuristic argument yields predictions for the limit $\psi$ in many ``locally
tree-like" models, as discussed in section \ref{sec-different}.  
However for the proofs in this paper we exploit special structure of our model, and
it seems technically challenging to find rigorous proofs in the broader settings of 
section \ref{sec-different}.

\subsection{The complete graph with random edge-lengths}
Our probability model for a random $n$-vertex network starts with the complete graph
and assigns independent Exponential(rate $1/n$)
random lengths $L_{ij} = L_{ji} = L_e$ to the
${n \choose 2}$ edges $e = (i,j)$.
This model (which we denote by $\GG_n$) and minor variants (uniform$(0,1)$ lengths;
complete bipartite graph) have been
studied in various
contexts, for instance the length of minimum spanning tree
\cite{fri85}, Steiner tree \cite{bol-steiner}, minimum matching
\cite{me94,LW04,MP87,NPS05} and traveling salesman
tour
\cite{MP86,me103,wast06a}. 
Note that our scaling convention $\Ex L_e = n$ makes lengths a factor $n$ larger
than in most of the earlier literature.
Most closely related to the present paper is the work of
Janson \cite{janson-123} and van der Hofstad et al. \cite{HHV02}  who studied
several aspects of the distances $D_n(i,j)$: see also W{\"a}stlund \cite{wast-spp} for connections with minimum matching.
In particular, it is known (\ref{Dn}) that 
$\Ex \bar{D}_n = (1 + o(1)) \log n$
so that we use 
$\bar{d}_n := \log n$ 
to scale edge-flows.

\subsection{The main result}
To fix notation, the vertex-set is
$[n] := \{1,2,\ldots,n\}$.
All quantities in the $n$-vertex model $\GG_n$ depend on $n$; our
notation makes $n$ explicit only where helpful. 
 For each ordered pair $(i,j)$ write $\bpi(i,j)$ for the shortest path
(considered here as a set of directed edges $e$) from $i$ to $j$. Define
 \begin{equation}
F_n(e) = \sfrac{1}{n} \sum_{i \in [n]} \sum_{j \in [n], j \neq i} \ind \{e \in
\bpi(i,j)\}
\label{Fndef}
\end{equation}
so that  $F_n(e)$ is the total flow across the directed edge $e$ in the specified
direction,
 when a flow of volume $1/n$ is put along the shortest path between each
ordered vertex pair.  
 Write $\#$ for
cardinality.
 \begin{Theorem}
\label{T1}
As $n \to \infty$ for fixed $z>0$,
\begin{equation}
\sfrac{1}{n} \# \{e: F_n(e) > z \log n\}
\to_{L^1}
G(z):=
  \int_0^\infty \Pr(W_1W_2e^{-u} > z) \ du
 \label{T1-L1}
\end{equation}
 where $W_1$ and $W_2$ are independent
 Exponential($1$).
 In particular
\begin{equation}
 \sfrac{1}{n} \Ex \# \{e: F_n(e) > z \log n\} \to
 G(z) . 
 \label{T1-Ex}
\end{equation}
In more detail, for $\bar{d}_n:= \log n$ the random empirical measure $\psi_n$ at (\ref{psi-n}) converges to
the non-random measure $\psi$ which is the ``distribution" of
$( U_\infty,W_1W_2 e^{-U_\infty})$ when $U_\infty$ is uniform on
$(0,\infty)$ and independent of $(W_1,W_2)$.
 \end{Theorem} 
In the final assertion we wrote ``distribution" because $\psi$ is a 
$\sigma$-finite distribution.  
As explained in section \ref{sec-complete-proof},  ``convergence of $\psi_n$" means $L^1$ convergence over the vague topology.  
The appearance of a $\sigma$-finite limit is not surprising, because edges
 of fixed large length carry a flow which is small but non-negligible
compared
 to flow across edges of length $1$.
 Note the anticipated identity (\ref{psi-1}) holds because 
 $\int_0^\infty ue^{-u} \ du = 1$.   
 Note also that the scaling of edge-lengths in $\GG_n$ does not affect the conclusions 
 (\ref{T1-L1},\ref{T1-Ex}) which remain true if edge-lengths have Exponential($1$) or 
 Uniform$(0,1)$ distribution.

See section \ref{sec-Gz} for further discussion of the function $G(z)$. In
particular its tail behavior is
 a stretched exponential
\begin{equation}
 G(z) = \exp( - z^{\frac{1}{2} + o(1)}) \mbox{ as } z \to \infty
 \label{Gzasy}
 \end{equation}
 rather than an ordinary exponential as one might have guessed.  
 Section \ref{sec-vertexflows} states the analog of Theorem \ref{T1} for the distribution of flows through {\em vertices} instead of edges.

\subsection{A heuristic argument}
\label{sec-heur}
Here is a heuristic argument for why the limit is this
particular function $G(z)$.
Consider a short edge $e$, that is an edge of length $O(1)$.  Suppose there are
$W_e^\prime(\tau)$ vertices within a fixed large distance $\tau$
of one end of $e$, and
$W_e^{\prime \prime}(\tau)$ vertices within distance $\tau$
of the other end.
A shortest-length path between
distant vertices which passes through $e$ must enter and exit the
region above via some pair of vertices in the sets above (see Figure 2), and there
are
$W_e^\prime(\tau) W_e^{\prime \prime}(\tau)$
such pairs.
The dependence on the length $L_e$ is more subtle.
By the Yule process approximation (Lemma \ref{LYule})
the number of vertices within distance $r$ of an initial vertex
grows as $e^r$, and it turns out that the flow through $e$ depends on $L_e$
as $\exp(-L_e)$ because of the availability of alternate possible shortest paths.
So flow through $e$ should be proportional to
$W_e^\prime(\tau) W_e^{\prime \prime}(\tau)
\exp(-L_e)$.
But (again by the Yule process
approximation, Lemma \ref{LYule})
for large $\tau$ we have
$e^{- \tau}W_e^\prime(\tau)$
has approximately an Exponential($1$) distribution $W_1$.
And as $n \to \infty$ the normalized distribution
$n^{-1} \#\{e: L_e \in \cdot\}$
over directed edges converges to the $\sigma$-finite distribution
of $U_\infty$.
This is heuristically how the limit joint distribution
$(U_\infty, W_1W_2 \exp(-U_\infty))$
arises.

Our proof of Theorem \ref{T1} is % (sections \ref{sec-proof} and \ref{sec-var}) is
essentially just a formalization of the heuristic argument using
explicit calculations exploiting the special structure of our
random network model.  
But we exploit a variety of tools to handle the details.  
For proving (section \ref{sec-proof}) the ``expectation" assertion (\ref{T1-Ex})
the key idea is
\begin{itemize}
\item analyzing the behavior of the percolation (flow from vertex $1$) process in a given neighborhood (sections \ref{sec-pt} and \ref{sec-condmean})
\end{itemize}
but we also use 
\begin{itemize}
\item the Yule process local approximation (section \ref{sec-Yule})
\item a martingale property (section \ref{sec-martingale}) 
\item a general weak law of large numbers for local functions on $\GG_n$ (section \ref{sec-WLLN}).
\end{itemize}
For proving (section \ref{sec-var})
the $L^1$ convergence assertion (\ref{T1-L1}),
we need to study the joint behavior of two shortest paths 
$\bpi(1,2), \bpi(3,4)$.  
This involves somewhat intricate conditioning arguments.  
The key ideas are
\begin{itemize}
\item finite-$n$ bounds for mean intensities of short paths (section \ref{joint-intensity})
\item 
the size-biased Yule process (section \ref{sec-sizebias})
\item conditional on existence of a given short path from vertex $1$,
the process of numbers of vertices within distance $t$ from vertex $1$
grows as a size-biased Yule process (section \ref{sec-sizebias}).
\end{itemize}

\section{Proofs}
\label{sec-proof}

\subsection{Preliminaries}
\label{sec:prelim}
Exponential$(\lambda)$ and Geometric$(p)$ denote the exponential and geometric
distributions in their usual parametrizations.

Here we collect without proof some standard properties of the
random network model $\GG_n$.
For fixed $n$ % a fixed vertex (vertex $1$, say)
and for $t \geq 0$ define
\begin{equation}
N_n(t):=
\mbox{ number of vertices within distance $t$ from vertex $1$}
\label{def-Nnt}
\end{equation}
where we include vertex $1$ itself;
\[ S_{n,k} := \min \{t: N_n(t) = k\}, \quad
 1 \leq k \leq n-1 \]
so that $S_{n,k+1}$ is the distance from vertex $1$ to the $k$'th nearest distinct vertex.
Then
\begin{equation}
\left(S_{n,k+1}-S_{n,k}, 1 \leq k \leq n-1\right) \mbox{ are
independent Exponential}(\sfrac{k(n-k)}{n})
 .  \label{Snk}
\end{equation}
Because the distance $D(1,2)$ is distributed as  $S_{n,V}$ for $V$ uniform
 on $\{2,3,\ldots,n\}$, it is straightforward (see e.g. \cite{MR0100305,HHV02}
for similar calculations)
to use (\ref{Snk}) to write exact formulas for the 
mean, variance and generating function of $D(1,2)$ 
and then deduce the $n \to \infty$ limit behavior
\begin{equation}
\Ex D(1,2) - \log n \to c_1, 
\quad \var D(1,2) \to c_2, 
\quad D(1,2) - \log n \cd D_\infty 
\label{D12-limits}
\end{equation}
for finite constants $c_1,c_2$ and a distribution $D_\infty$ 
discussed further in section \ref{sec-k-paths}.  
Note that the average vertex-vertex distance $\bar{D}_n$ at (\ref{feD}) has the same mean, but not the same distribution, as $D(1,2)$, so
\begin{equation}
\Ex \bar{D}_n = \log n + O(1)
\mbox{ as } n \to \infty .
\label{Dn}
\end{equation}

There is a natural mental picture of
(first passage) percolation, in which at time $0$ there is water
at vertex $1$ only, and the water spreads along edges at speed $1$.
So $N_n(t)$ vertices have been wetted by time $t$.
Each vertex $j \neq 1$ is first wetted via some edge $(i(j),j)$,
and this collection of directed edges forms the
{\em percolation tree} rooted at vertex $1$.
The ``flow" in Theorem \ref{T1} from vertex $1$ goes along the
edges of this percolation tree, and no other edges.

Associated with the percolation process is a filtration $(\FF_t)$,
where $\FF_t$ is the information known at time $t$, illustrated
informally as follows. Write $T_i$ for the wetting time of vertex
$i$. Then the values of $T_i$ for which $T_i \leq t$ are in
$\FF_t$. On the event $\{T_i < t < T_j\}$, the information in
$\FF_t$ about $L_{ij}$ is that $L_{ij} > t-T_i$. By the memoryless
property of the exponential distribution, on the event above the conditional
distribution of $L_{ij} - (t-T_i)$ given $\FF_t$ is
Exponential$(1/n)$,
and the $(L_{ij})$ are conditionally independent given $\FF_t$.
More elaborate versions of this memoryless property 
appear in Lemmas \ref{Lneigh} and \ref{LNAF}.

An obvious consequence of the Exponential$(1/n)$ distribution of edge-lengths $L_{ij}$
is that 
\[ n \Pr (L_{ij} \leq \tau) \to \tau \mbox{ as } n  \to \infty . \]  
In words, this says that the measure 
$n \Pr(L_{ij} \in \cdot)$ converges {\em vaguely} to Lebesgue measure on $(0,\infty)$.

\subsection{A martingale property}
\label{sec-martingale}
In the percolation process above,
write $\WW_n(t)$ for the \underline{set} of vertices wetted by time $t$.
The following martingale property turns out to be useful.
\begin{Lemma}
\label{Lmartin}
Let $W_n$ be a stopping time for the percolation process on $\GG_n$.
For each $v \in \WW_n(W_n)$ let $Y(v)$ be the number of vertices
$j \in [n]$ such that, in the shortest path from $1$ to $j$,
the last-visited vertex of $\WW_n(W_n)$ is vertex $v$.
Then
\[ \sfrac{1}{n} \Ex ( Y(v) | \FF_{W_n}) = \sfrac{1}{N_n(W_n)} . \]
\end{Lemma}
\proof
Define $Y(v,t)$ as $Y(v)$ but counting only vertices $j$ which are wetted by time $t$.
As $t$ increases, whenever a new vertex $j$ is wetted via some edge $(i(j),j)$,
the predecessor vertex $i(j)$ is a uniform random element of $\WW(t-)$.
It follows easily that the process
$\frac{Y(v,t)}{N_n(t)}, \ t \geq W_n$ 
is a martingale.
The optional sampling theorem shows
\[ \Ex \left( \left. \frac{Y(v,\infty)}{N_n(\infty)} \right| \FF_{W_n} \right) 
= \frac{Y(v,W_n)}{N_n(W_n)} = \frac{1}{N_n(W_n)} . \]
But 
$\frac{Y(v,\infty)}{N_n(\infty)} = \frac{Y(v)}{n}$.
\qed

\subsection{The Yule process approximation}
\label{sec-Yule}
The Yule process
$(N_\infty(t), \ 0 \leq t < \infty)$
is the population at time $t$ in the continuous-time branching process
started with one individual, in which individuals live forever
and produce offspring at the times of a Poisson($1$) process.
Writing
\[ S_{\infty,k} := \min \{t: N_\infty(t) = k\} \]
we have
\begin{equation}
\left(S_{\infty,k+1}-S_{\infty,k}, 1 \leq k \leq \infty \right)
\mbox{ are independent Exponential}(k)%(\sfrac{k(n-k}{n})
\label{Sik}
\end{equation}
and so the Yule process is the natural
$n \to \infty$ limit of
the process (\ref{Snk}) associated with the percolation tree.
We quote some standard facts about the Yule process.
\begin{Lemma}
\label{LYule}
(a) $N_\infty(t)$ has Geometric($e^{-t}$) distribution.\\
(b) $e^{-t}N_\infty(t)$ is a martingale which is bounded in $L^2$,
and $e^{-t} N_\infty(t) \to W$ a.s. and in $L^2$ as $t \to \infty$,
where $W$ has Exponential($1$) distribution.
\end{Lemma}

It is intuitively clear that that the local structure of $\GG_n$ relative to one vertex 
converges to the Yule process.  
Abstractly \cite{me101}, we call this notion of convergence of random graphs 
{\em local weak convergence} and the limit structure 
(the Yule process regarded as a ``spatial" graph) is called the PWIT.  
But rather than work abstractly we will state only the more concrete  consequences needed, 
such as the next lemma.
\begin{Lemma}
\label{Lkt}
Fix $k \geq 1$ and $t < \infty$. 
For $1 \leq i \leq k$ let $N^{(i)}_n(t)$ be the number of vertices of $\GG_n$ within distance $t$ from vertex $i$.  Then as $n \to \infty$
\[ (N_n^{(1)}(t),\ldots,N_n^{(k)}(t)) \cd (N_\infty^{(1)}(t),\ldots,N_\infty^{(k)}(t)) \]
where the limits $N_\infty^{(i)}(t)$ are independent Geometric$(e^{-t})$.
\end{Lemma}
\proof 
For $k = 1$ this follows from (\ref{Snk},\ref{Sik}) and Lemma \ref{LYule}(a). 
For general $k$, use the natural conditioning argument.
\qed

The following technical lemma shows one way in which the ``exponential growth
with rate $1$" property of the Yule process (Lemma \ref{LYule}(b)) translates to the
percolation process.
\begin{Lemma}
\label{LNs}
Let $W_n$ be a randomized stopping time for the percolation
process on $\GG_n$.
Fix $\eps > 0, \sigma < \infty$
and a sequence 
$(\omega_n)$ such that
$\omega_n \to \infty$ with 
$\omega_n \leq n^{1/2}$.
Then
as $n \to \infty$
\[ \Pr \left( \left. 1 - \eps
\leq \frac{N_n(W_n+\sigma)}{e^\sigma N_n(W_n)} \leq 1+ \eps
\right| \FF_{W_n} \right) \to 1 \] uniformly on $\{\omega_n \leq
N_n(W_n) \leq n/\omega_n\}$.
\end{Lemma}
\proof
It is enough to show this holds conditionally on $N(W_n)$,
that is to show 
\begin{equation}
 \Pr \left( \left. 1 - \eps
\leq \frac{N_n(t_n+\sigma)}{e^\sigma k_n} \leq 1+ \eps
\right| W_n=t_n, N_n(W_n) = k_n \right) \to 1 \label{Nws} % uniformly on $\{\omega_n \leq
\end{equation}
whenever $k_n \to \infty, % k_n \leq n^{1/2}$.
n/k_n \to \infty$.
Note that the value of $t_n$ does not affect the conditional probability.

First note that by (\ref{Snk},\ref{Sik}) we can couple 
$N_n(\cdot)$ and $N_\infty(\cdot)$
by constructing each from the same i.i.d. Exponential($1$) sequence
$(Y_i)$
via
\begin{eqnarray}
 S_{\infty,k} &=&
\sum_{i=1}^{k-1}
\sfrac{1}{i} Y_i; \quad 
N_\infty(t) = \max \{k: S_{\infty,k} \leq t \} \nonumber \\
 S_{n,k} &=& 
\sum_{i=1}^{k-1}
\sfrac{n}{i(n-i)} Y_i; \quad 
N_n(t) = \max \{k: S_{n,k} \leq t \} . \label{SScouple}
\end{eqnarray}
>From this coupling we see that for
$k_n < m_n$ with
$k_n \to \infty, \ n/m_n \to \infty$
we have
\begin{equation}
 \frac
{S_{n,m_n} - S_{n,k_n}}
{S_{\infty,m_n} - S_{\infty,k_n}}
\to 1 \mbox{ a.s.} \label{SS}
\end{equation}
By the homogeneous branching property of the Yule process
and Lemma \ref{LYule}(a),
conditional on $\{N_\infty(t_n) = k_n\}$,
we can represent 
$N_\infty(t_n + \sigma)$
as the sum of $k_n$ independent Geometric($e^{-\sigma}$) r.v.'s,
and so (still conditionally)
$
\frac{
N_\infty(w_n + \sigma)
}{k_n e^\sigma} \to 1$
in probability.
In terms of 
$(S_{\infty,k})$
this is equivalent to the 
(now unconditional)
property that
\[ \mbox{ if } k_n \to \infty, \ 
m_n \sim k_n e^\sigma 
\mbox{ then } 
S_{\infty,m_n} - S_{\infty,k_n} \cp \sigma . \]
Now by (\ref{SS}) and assumptions on $(k_n)$ we see
\[ \mbox{ if } k_n \to \infty, \ 
m_n \sim k_n e^\sigma 
\mbox{ then } 
S_{n,m_n} - S_{n,k_n} \cp \sigma . \]
This holds for each fixed $\sigma > 0$.
Translating this back into an assertion about
$(N_n(\cdot))$
establishes (\ref{Nws}).  
\qed

Recall that $D(i,j)$ denotes vertex-vertex distances in $\GG_n$.
\begin{Lemma}
\label{LBCD}
For disjoint subsets $B,C \subset [n]$ and for any $d > 0$,
\[ \Ex \# \{(i,j) \in B \times C: \ D(i,j) \leq d\}
\leq \#B \#C e^d/(n-1) . \]
\end{Lemma}
\proof
By linearity and symmetry we reduce to the case where $B$ and $C$
are singletons.
Then the left side equals
$\Pr(D(1,2) \leq d) = (n-1)^{-1} (\Ex N_n(d) - 1)$.
By the coupling (\ref{SScouple}) to the Yule process, and Lemma \ref{LYule}(a),
$ \Ex N_n(d) \leq \Ex N_\infty(d) = e^d $. 
\qed

\subsection{Local structure in the $n$-vertex model}
\label{sec-Ls}
In this section we give a result (Lemma \ref{Lneigh})
describing the global structure of $\GG_n$ conditional
on a given local structure.
The actual result is obvious once stated, but requires some notational effort
to set up.

Fix a real $\tau > 0$.
Let $\bt$ be a finite unlabelled tree with edge lengths, with the following properties (see Figure 1).\\
(i) There is a distinguished directed edge, whose % which serves to label
end vertices can then be labelled as $(v_L,v_R)$, defining a partition of all the vertices
of $\bt$ as $V(\bt) = L(\bt) \cup R(\bt)$.
Here $L$ and $R$ are mnemonics for {\em left} and {\em right}, and the distinguished edge is directed left-to-right.\\
(ii) Every vertex in $L(\bt)$ is within distance $\tau$ from $v_L$,
and every vertex in $R(\bt)$ is within distance $\tau$ from $v_R$.

\vspace{0.2in}
\begin{figure}
\begin{center}

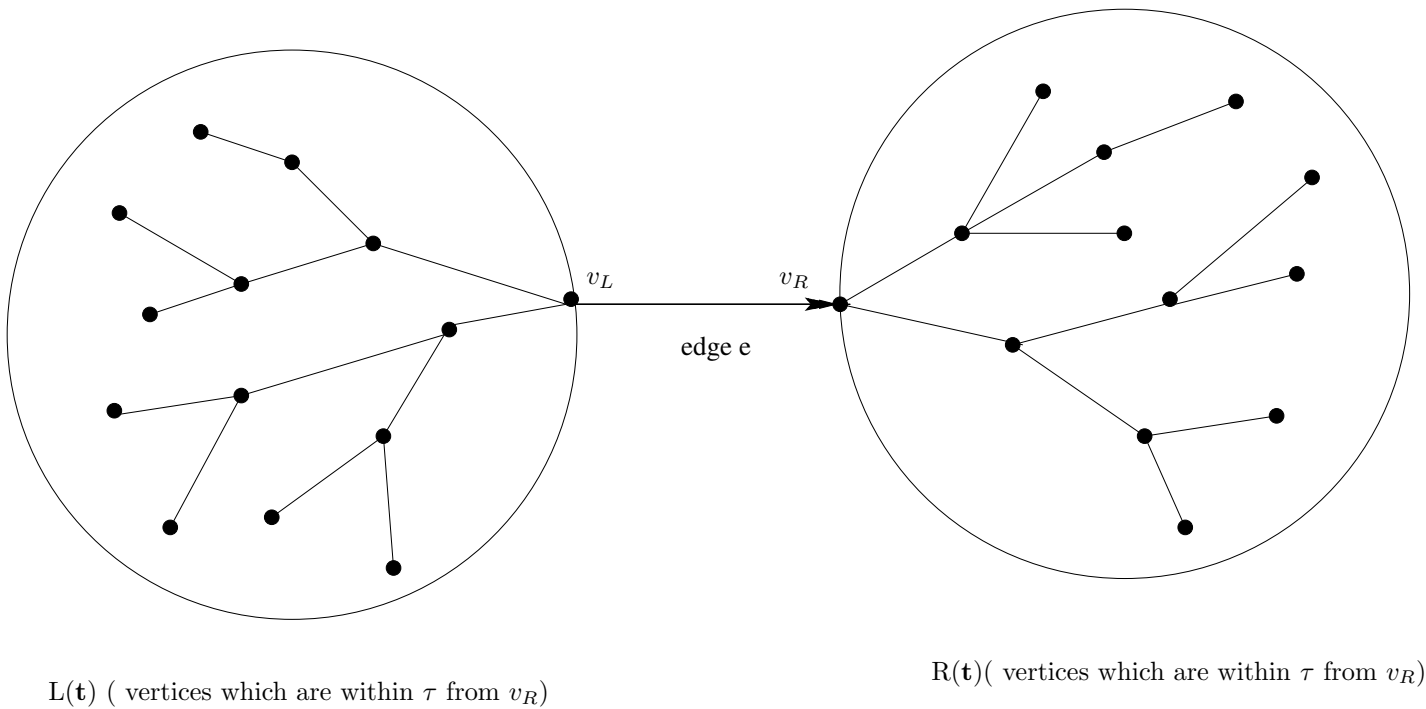

% \vspace{-0.5in}

\caption{\bf Neighborhood about an edge} \label{nbhd about an
edge}
\end{center}
\end{figure}

\vspace{0.2in}
\noindent
Write $\bT_\tau$ for the set of such trees $\bt$.
For such $\bt$, write
$\ell(e)$
for the length of the distinguished edge $e$.
And for each vertex $v \in L(\bt)$ write
$b(v) = \tau - D(v,v_L)$ as the
``distance to boundary" 
(as in Figure 1, we envisage a boundary drawn at distance  at distance $\tau$ from $v_L$ and from $v_R$).
Similarly 
for $v \in R(\bt)$ write
$ b(v) = \tau - D(v,v_R)$.  %\quad v \in R(\bt) . \]
Finally, given $\bt$ and given a subset $A \subseteq [n]$
with $\# A = \# V(\bt)$,
let $\bt_A$ denote some labelling of the vertices of $\bt$ by distinct
labels from $A$.

Now consider the random network $\GG_n$.
We occasionally want to regard an edge $(i,j)$ of $\GG_n$
as a point-set, so that a point on the edge is at some distance
$0 < u < L_{ij}$ from vertex $i$ and at distance $L_{ij} - u$
from vertex $j$, with this notion of distance extending in the natural way
to distances between a point on an edge and a distant vertex.
Fix $\tau$ and vertices $v_L,v_R \in [n]$.
Define the ``neighborhood"
$\NN_\tau(v_L,v_R)$ as the subgraph of $\GG_n$
whose vertex-set consists of
vertices $v$ for which
$\min(D(v,v_L), D(v,v_R)) \leq \tau$.
Its edge-set is the subset of edges of $\GG_n$ such that,
for every point along the edge, the distance to the closer of
$\{v_L,v_R\}$ is at most $\tau$.
It also contains by fiat the distinguished directed edge $(v_L,v_R)$.

In general
$\NN_\tau(v_L,v_R)$
need not be a tree, but if it is a tree then clearly
it is a tree of the form $\bt_A$ for some $\bt \in \bT_\tau$
and some $A \subseteq [n]$.
In this case we can define $b(i)$ (meaning distance to boundary of neighborhood) as above for 
vertices $i$ of 
$\NN_\tau(v_L,v_R)$,
and we define $b(i) = 0$ for other vertices of $\GG_n$.
\begin{Lemma}
\label{Lneigh}
Fix $n, \tau, v_L, v_R$ and $\bt_A$.
Then conditional on
$\NN_\tau(v_L,v_R)
= \bt_A$,
the lengths $(L_{ij})$ of the edges of $\GG_n$ which are not
edges of
$\NN_\tau(v_L,v_R)$
are independent r.v.'s % and for these edges the r.v. 
for which
$L_{ij} - b(i) - b(j)$ has Exponential($1/n$) distribution.
\end{Lemma}
\proof
Saying that $(i,j)$ is not an edge of 
$\NN_\tau(v_L,v_R)$
is saying that 
$L_{ij} > b(i) + b(j)$.
The edge-lengths $(L_{ij})$ are a priori independent Exponential($1/n$), and 
conditioning on all these inequalities leaves them independent with the
stated distributions.

\subsection{The percolation tree on a neighborhood}
\label{sec-pt}
We now come to the central idea of the proof, which is to
study how the percolation tree behaves on a given neighborhood.
Until further notice we adopt the setting of Lemma \ref{Lneigh} and work conditionally on
$\NN_\tau(v_L,v_R)
= \bt_A$.
 Let $\conde$ and $ \condp$ denote respectively the conditional expectation and
conditional probability operations;
and we will use tilde notation $\wL_{ij}$ to denote
{\em unconditioned} quantities.
Thus the conclusion of Lemma \ref{Lneigh} can be rewritten as follows.
Starting with independent Exponential($1/n$) r.v.'s $(\wL_{ij})$ 
we can construct the conditioned lengths $(L_{ij})$ as 
\begin{equation}
L_{ij} = \wL_{ij} + b(i) + b(j), \quad
(i,j) \notin  \NN_\tau(v_L,v_R) . \label{LL}
\end{equation}  
This provides a coupling of the unconditioned and conditioned lengths.

\vspace{0.2in}
\begin{figure}
\begin{center}

 \centerline{\resizebox{.8\linewidth}{!}{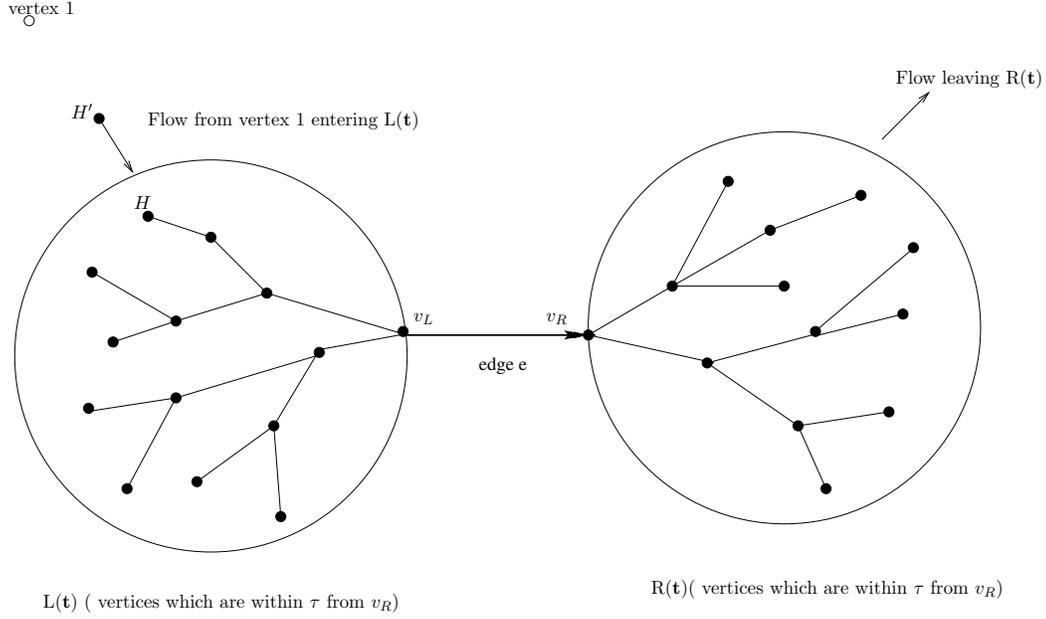}}

\caption{{\bf Flow from 1 passing through the neighborhood.}  Flow enters at time $T_L$ along some edge 
$(H^\prime,H)$, and exits at time $T_L + 2\tau + \ell(e)$.
} \label{flow entering}
\end{center}
\end{figure}

Now consider the percolation process started at vertex $1$, and assume
$1 \notin A$.
Let $T_L$ be the first time
(in the conditioned model)
that the percolating water gets to some point on an edge at distance $\tau$
from $v_L$
without passing along the distinguished edge.
Define $T_R$ similarly and then set
$T = \min(T_L,T_R)$.
So at time $T$ the water is at distance $b(H)$ from some random
vertex $H$ of $A$.   
See Figure \ref{flow entering}.
\begin{Lemma}
\label{Lexact}
(a)
The (conditioned) distribution of $T$ is the same as the
(unconditioned) distribution of $\wT_A$, the first time in
percolation on $\GG_n$ that some vertex in $A$ is wetted.
And the random vertex $H$ is distributed uniformly on $A$,
independent of $T$.
\\(b)
$N(T)$ is distributed as the smallest of $\#A$ uniform random samples
without replacement from $\{2,3,\ldots,n\}$.
\\(c)
Let $T^\prime$ denote the second time that the percolating water gets
to within distance $\tau$ of either $v_L$ or $v_R$ 
along some path which has not previously 
hit $\bt_A$.
Then $(N(T),N(T^\prime))$ has the joint distribution of the
smallest and the second smallest of $\#A$ uniform random samples
without replacement from $\{2,3,\ldots,n\}$.
\end{Lemma}

{\em Proof.} Use (\ref{LL}) to construct the conditioned process
from the unconditioned process. In the unconditioned process it is
clear by symmetry that $H$, the first vertex of $A$ wetted, is
uniform on $A$ and independent of $\wT_A$. Obviously $H$ is
reached along some edge $(H^\prime,H)$ with $H^\prime \notin A$.
In the conditioned process, at time $\wT_A$ the percolating water
has reached distance $b(H) = b(H) + b(H^\prime)$ from $H$, and
hence is distance $\tau$ from either $v_L$ or $v_R$ (whichever is
closer to $H$). So $T = \wT_A$ in the coupling. This gives (a). Parts (b) and (c)
are similar. 
\qed

The next lemma formalizes the idea 
``what do we know about edge-lengths at time $T_L$?"
As described above, on
$\{T_L<T_R\}$ 
there is some vertex $H \in L(\bt_A)$
such that $T_H - T_L = b(H)$, and vertex 
$H$ gets wetted via some edge $(H^\prime,H)$
where $H^\prime$ is in the set 
$\WW(T_L)$
of vertices wetted by time $T_L$.
Arguing as in Lemma \ref{Lneigh} shows 
\begin{Lemma}
\label{LNAF}
Conditional on 
$\NN_\tau(v_L,v_R)
= \bt_A$
and conditional on 
$\sigma(\FF_{T_L},H^\prime,H)$,
on the event 
$\{T_L<T_R\}$,
the collection of edge-lengths $(L_{ij})$
as $(i,j)$ runs over all edges except\\
(i) $(i,j)$ an edge of 
$\NN_\tau(v_L,v_R)$
\\ (ii) 
$i,j \in \WW(T_L)$\\ (iii)
$(i,j) = (H^\prime, H)$\\
are independent with distributions\\
\[ L_{ij} = \wL_{ij} + b(i) + b(j), \quad
(i,j) \notin  % \label{LL}
\WW(T_L) 
\]
\[ 
L_{ij} = \wL_{ij}
+ (T_L - T_i) + b(j), \quad 
i \in \WW(T_L) . \]
\end{Lemma}

Recall we are conditioning on 
$\NN_\tau(v_L,v_R)
= \bt_A$. 
Write
\[ \sigma := 2 \tau + \ell(e) \]
where $e$ is  the distinguished edge 
$e = (v_L,v_R)$ of $\bt_A$.
Estimating the mean flow through $e$
is tantamount to estimating 
the mean of the random variable 
\begin{equation}
  M_1 := 
\# \{2 \leq j \leq n: \ e \in \bpi(1,j) \} \label{M1-def}
\end{equation}
counting the  number of vertices
$j \in [n]$ with $j \neq 1$
such that the shortest path from $1$ to $j$ passes through
$e$.
To analyze $M_1$, we consider 
the set $\RR^*(\bt_A)$ of vertices in $R(\bt_A)$ which are first reached via edge $e$: 
\[ \RR^*(\bt_A) := \{v \in R(\bt_A): e \in \bpi(1,v)\} . \] 
Note that from the definitions of $T_L$ and $T_R$ 
\begin{eqnarray*}
\mbox{ if } T_R  \leq T_L & \mbox{ then } & \# \RR^*(\bt_A) = 0  \mbox{ and } M_1 = 0\\
\mbox{ if } T_L < T_R \leq  T_L + \sigma & \mbox{ then }  & 0 \leq \# \RR^*(\bt_A) \leq \# R(\bt_A) \\
\mbox{ if } T_R > T_L + \sigma  & \mbox{ then } &  \# \RR^*(\bt_A) =  \# R(\bt_A) .
\end{eqnarray*}
\begin{Lemma}
\label{lem:expec-0}
$\mathbb{E}_\tau(M_1|\FF_{T_L+\sigma}) =\frac{n \ \# \RR^*(\mathbf{t}_A)}{N(T_L + \sigma)}$ .
\end{Lemma}
\proof
In the notation of Lemma \ref{Lmartin}
\begin{equation}
 M_1 = \sum_{v \in \RR^*(\bt_A)} Y(v)  .  \label{M1Y}
\end{equation}
Applying Lemma \ref{Lmartin} with $W = T_L + \sigma$
(the fact we are working conditional on $\NN_\tau(v_L,v_R)
= \bt_A$ doesn't affect the martingale property after time $W$) 
gives the second equality below:
\[ \sfrac{1}{n} \conde (M_1|\FF_W)
=  \sum_{v \in \RR^*(\bt_A)} \sfrac{1}{n} \Ex_\tau(Y(v) | \FF_W) 
= \sfrac{\# \RR^*(\bt_A)}{N(W)}
. \]
\qed

\subsection{The conditioned mean flow}
\label{sec-condmean}
We now start studying $n \rightarrow \infty$ asymptotics.  Recall the definition (\ref{Fndef})  of the normalized flow $F_n(e)$ across an edge $e$ of $\GG_n$.  The next result calculates the expected flow conditional on the neighborhood structure of $\GG_n$ around $e$.
  \begin{Proposition}
 \label{Pcond-mean}
Fix $\tau$ and $\mathbf{t} $ $\in $ $\mathbf{T}_\tau$ and write $\sigma = 2\tau + \ell(e)$ where $e$ is the distinguished edge $e$ of $\bt$. Let $v_L^n \neq v_R^n  $  $\in $ $[n]$ and let $\{v_L^n, v_R^n\} \subseteq A_n \subset [n]$ satisfy $\#A_n = \#\mathbf{t}$.   Then as $n \rightarrow \infty$, setting $e_n = (v_L^n , v_R^n )$,
\[\expec_\tau(F_n(e_n)) = (1+o(1)) \#L(\bt) \#R(\bt) e^{-\sigma} \log{n} . \]
 \end{Proposition}
Here $\expec_\tau(\cdot)$ denotes $\expec(\mbox{ $ \cdot  $} | \NN_\tau(v_L^n ,v_R^n ) = \bt_{A_n})$.  
All quantities except $\bt, \tau, \sigma$ depend on $n$, though the dependence is often not made explicit in notation 
in the proof below.

{\bf Proof.}  Assume vertex $1 \notin A$. Recall the definition of $M_1$  from (\ref{M1-def}). 
We shall show that 
\begin{equation}
\label{eqn:m}
\expec_\tau(M_1) = (1+o(1)) \#L(\bt) \#R(\bt) e^{-\sigma} \log{n} .
\end{equation}
 In terms of flows of volume $1/n$ between each vertex-pair, the (conditional) mean contribution to the flow $F_n(e)$ through the distinguished edge arising from flow started at vertex $1 \notin A$ equals $n^{-1} \Ex(M_1)$.  The same contribution arises from each of the $n- \# \bt$ possible starting  vertices $v \notin A$. For $v \in A$ the flow through e is trivially bounded by 1. So to prove Proposition \ref{Pcond-mean} it is enough to prove (\ref{eqn:m}).
 
 Fix a sequence $\omega_n \rightarrow \infty$, with $\omega_n = o(\log{n})$. 
  Start the first passage percolation process from vertex 1. 
  Recall that $T_L$ denotes the first time the flow is 
  within distance $\tau$  from $v_L$, 
and that by time $T_L + \sigma$ the flow has wetted every vertex in $\bt_A$.
  We shall show that the % contributions to $\conde(M_1)$ when the flow from 1 
   dominant contribution to $\conde(M_1)$ is from  the ``good" event 
  \begin{equation}
   G^* := \{\omega_n\leq N_n(T_L) \leq n/\omega_n\} \cap \{T_R > T_L + \sigma\}. \label{G*def}
   \end{equation}

For the details, 
observe that we can apply Lemma \ref{LNs} to the percolation process on edges excluding the edges of $\bt_A$ and deduce
\begin{equation}
\label{eqn:exponential}
\condp \left( \left. 1-\eps \leq \frac{N_{\mbox{{\tiny avoid}}}(T_L + \sigma)}{e^\sigma N(T_L)} \leq 1+\eps \right|  \FF_{T_L}\right) \rightarrow 1
\mbox{ uniformly on $\{\omega_n \leq N(T_L) \leq n/\omega_n\}$} 
\end{equation}
where $N_{\mbox{{\tiny avoid}}}(t)$ is the number of vertices wetted by time $t$ via paths which use no edge of $\bt_A$. 
For the rest of the argument we work on the event 
$\{T_L < T_R\}$ 
(otherwise, $M_1 = 0$). 
Recall (Figure 2) that flow enters the neighborhood at time $T_L$ along some edge 
$(H^\prime,H)$.
We claim: conditional on $\sigma(\mathcal{F}_{T_L} , H ,H^\prime)$, on the event $\{T_L < T_R\}$, the expected number of vertices of $R(\bt_A)$ wetted before time $T_L + \sigma$ by paths \underline{not} using the distinguished edge is at most
\begin{equation}
\label{bad-event} 
(N(T_L) + \#L(\bt)) \#R(\bt) e^\sigma / (n-1) .
\end{equation}

This follows from Lemmas \ref{LNAF} and \ref{LBCD}  applied to $\mathcal{W}(T_L) \cup L(\bt_A)$ and $R(\bt_A)$, because the former lemma implies that the conditioned edge lengths can only be longer than the unconditioned edge-lengths in the latter lemma. 
Note that the expectation (\ref{bad-event}) tends to 0 on $\{N(T_L) \leq n/\omega_n \}$, so that
\begin{equation}
\label{eqn:before-sigma}
\condp(T_R >T_L + \sigma|\mathcal{F}_{T_L}) \rightarrow 1 
\mbox{ uniformly on } \{ N(T_L) \leq n/\omega_n\} \cap \{T_L < T_R\} .
\end{equation}

Next let $N_{\mbox{{\tiny via}}}(t)$ denote the number of vertices outside $\bt_A$ which have been wetted by time $T_L+\sigma$ using some path via $\bt_A$.  Again using Lemmas \ref{LNAF} and  \ref{LBCD}, applied now to $L(\bt_A)$ and $[n] \setminus  A$ , we find
\begin{equation}
\conde(\ind\{T_L < T_R\} N_{\mbox{{\tiny via}}}(T_L+\sigma)| \FF_{T_L} ) \leq \#L(\bt) e^{\sigma} . 
\label{eq:via}
\end{equation}
>From the definitions, 
\[ 0 \leq N(T_L + \sigma) - N_{\mbox{{\tiny avoid}}}(T_L+\sigma) 
\leq N_{\mbox{{\tiny via}}}(T_L+\sigma) + \# \bt . \]
Combining this with (\ref{eqn:exponential}, \ref{eqn:before-sigma}, \ref{eq:via}) we deduce
\begin{equation}
\label{eqn:G-and-exponential}
\condp\left( \left. T_L + \sigma < T_R;1-\eps \leq \frac{N(T_L+\sigma)}{e^\sigma N(T_L)}\leq 1+\eps \right|\FF_{T_L} \right) \rightarrow 1
\mbox{ uniformly on $\{\omega_n \leq N(T_L) \leq n/\omega_n\}\cap \{T_L < T_R \}$} .
\end{equation}

Recalling the definition (\ref{G*def}) of $G^*$,   Lemma \ref{lem:expec-0} implies
\begin{equation}
\frac{\conde[M_1 \ind(G^*)]}{n \  \#R(\bt)} = \conde\left[ \frac{\ind(G^*)}{N(T_L + \sigma)} \right] . 
\label{eq:29}
\end{equation}

To obtain asymptotics for the right side, use the Lemma \ref{Lexact}(b) description of the distribution of $N(T_L)$ to conclude that 
\begin{equation}
\label{eqn:dist-of-N}
\condp\left(T_L < T_R , N_n(T_L) = m \right) = (1+o(1)) \frac{\#L(\bt)}{n-1}
\mbox{ uniformly on $\{2\leq m \leq n/\omega_n \}$ }.
\end{equation}
The harmonic sum estimate
$\sum_{\omega_n}^{n/\omega_n} j^{-1} = (1+o(1))
\log n$
leads to 
\[
\conde\left[\frac{1}{e^\sigma N(T_L)  } \ind\{\omega_n \leq N(T_L) \leq n/\omega_n\} \ind \{T_L + \sigma < T_R\}\right] = (1+o(1))\#L(\bt)e^{-\sigma}n^{-1}\log{n} 
 . \]
  Combine this with (\ref{eqn:G-and-exponential})  to get 
\[\conde\left[\frac{\ind(G^*)}{N(T_L+\sigma)} \right] =(1+o(1))\#L(\bt)e^{-\sigma}n^{-1}\log{n} \]
and thus by (\ref{eq:29}) 
\[\conde(M_1 \ind(G^*)) = (1+o(1)) \#L(\bt) \#R(\bt) e^{-\sigma} \log{n} . \]
Recalling that $M_1 = 0$ on $\{T_R < T_L\}$, we can write 
\[ \Ex_\tau [M_1]
= \Ex_\tau [M_1 \ind (G^*)]
+ \Ex_\tau [M_1 \ind (B_1)]
+ \Ex_\tau [M_1 \ind (B_2)]
+ \Ex_\tau [M_1 \ind (B_3)]  \]
for the ``bad" events
  \[B_1 :=  \{T_L <T_R\}  \cap \{N(T_L) > n/\omega_n \}\]
 \[B_2 := \{T_L <T_R\}  \cap \{N(T_L) < \omega_n \} \]
 \[B_3 :=  \{T_L <T_R< T_L +\sigma\} \cap \{\omega_n \leq N(T_L) \leq n/\omega_n \}  \]
and we need to check for each $B$ that 
$ \Ex_\tau [M_1 \ind (B)]
= o(\log n)$.
In each case we start by using Lemma \ref{lem:expec-0}.
\begin{eqnarray*}
\Ex_\tau [M_1 \ind (B_1)] &=& n \conde\left( \frac{\# \RR^*(\bt_A)}{N(T_L+ \sigma)} \ind\{T_L < T_R \} \ind\{ N(T_L) >  n/\omega_n\}\right)\\
&\leq& \#R(\bt) \omega_n  =  o(\log{n})
\end{eqnarray*}
because we chose $\omega_n = o(\log{n}) $.
\begin{eqnarray*}
\Ex_\tau [M_1 \ind (B_2)] &=& n \conde\left( \frac{\# \RR^*(\bt_A)}{N(T_L+ \sigma)} \ind\{T_L < T_R \} \ind\{ N(T_L) \leq \omega_n\}\right)\\
&\leq & n\#R(\bt) \condp\left(T_L < T_R , N_n(T_L)\leq \omega_n\right) 
\mbox{ by (\ref{eqn:dist-of-N})} \\
&=& (1+o(1))  \#R(\bt) \#L(\bt) \omega_n \\
&=& o(\log{n}) .
\end{eqnarray*}
For the third event,
\[
\Ex_\tau [M_1 \ind (B_3)] = n \conde\left( \frac{\# \RR^*(\bt_A)}{N(T_L+ \sigma)}
\ind \{T_L < T_R < T_L + \sigma\}
\ind\{\omega_n \leq N_n(T_L) \leq n/\omega_n\}\right)  . \]
By (\ref{bad-event}) we have
\[
\condp\left(\left. T_R <T_L + \sigma \right|\FF_{T_L}\right) \leq (N(T_L) + \#L(\bt)) \#R(\bt) e^\sigma / (n-1)
 .
\]
Because $N(T_L) \leq N(T_L+\sigma)$ and $\# \RR^*(\bt_A) \leq \# R(\bt)$,
writing $C=\{T_L < T_R\} \cap\{ \omega_n \leq N(T_L) \leq n/\omega_n\} $ we have 
\begin{eqnarray*}
 \Ex_\tau [M_1 \ind (B_3)] &\leq& n  \#R(\bt)  \conde \left( \frac{ (N(T_L) + \#L(\bt)) \#R(\bt) e^\sigma }{ n-1} . \frac{1}{N(T_L)} \ind(C) \right) \\
&\leq& (1+o(1))(\#R(\bt)+1)^2 \#L(\bt) e^\sigma .
\end{eqnarray*}
This completes the proof.
\qed

We record a minor rephrasing of Proposition \ref{Pcond-mean}.
\begin{Corollary}
\label{epi}
In the setting of Proposition \ref{Pcond-mean}, 
suppose vertices $1,2, \not\in A_n$.
Then 
\[
\condp(e_n \in \bpi(1,2)) =  (1+o(1)) \#L(\bt)\#R(\bt) \left[ \frac{\log{n}}{n} e^{-\sigma}\right]
 . \]
\end{Corollary}
\proof 
By symmetry over vertices $ j \not\in A_n \cup \{1\}$ we have 
$ \condp(e_n \in \bpi(1,2)) 
= \frac{\conde M_1^*}{n - 1 - \#\bt} $
where $M_1^*$ is defined as $M_1$ but excluding vertices $j \in A_n$.
Since $M_1 - M_1^* \leq \#\bt$, 
the corollary follows from (\ref{eqn:m}).
\qed

\subsection{Conditional variance of the flow}
\label{sec-cv}
In the setting of Proposition \ref{Pcond-mean} we want to show that the flow $F_n(e)$ is close to its conditional expectation, and the natural way to express this is via the conditional variance.

\begin{Proposition}
\label{Pcond-var} 
%Fix $B < \infty$
In the setting of Proposition \ref{Pcond-mean}.
\[ \frac{ \var_{\!\! \tau}  (F_n(e_n) )| \NN_\tau(v^n_L,v^n_R) = \bt_{A_n})}
{[\Ex_\tau  (F_n(e_n)| \NN_\tau(v^n_L,v^n_R) = \bt_{A_n})]^2} 
\leq \frac{1}{\# L(\bt)} +  \frac{1}{\# R(\bt)} + \frac{1}{\# L(\bt) \# R(\bt)} + o(1) 
\mbox{ as } n \to \infty . \]
\end{Proposition} 
This formulation emphasizes that the relative variance of the conditional distribution gets smaller as the  size of the neighborhood gets bigger.
We remark that Proposition \ref{Pcond-mean} alone 
(i.e. without Proposition \ref{Pcond-var})
is enough to prove the ``expectation" asssertion 
(\ref{T1-Ex}) of Theorem \ref{T1}.
Proposition \ref{Pcond-var} is needed for the $L^1$-convergence assertion (\ref{T1-L1}).

The key step in proving Proposition \ref{Pcond-var} is the following Proposition.
To set up notation, we may assume the edge $e_n$ is $(n-1,n)$ and the label set $A_n$ is  
$\{n- \# V(\bt) +1,\ldots n\}$.
Recall 
$\bpi(1,2)$ denotes the shortest path from $1$ to $2$. 
If this path uses $e_n$ then there is some entrance-exit pair 
$(\alpha,\beta) \in L(\bt) \times R(\bt)$ 
recording the first and last vertices of the neighborhood visited by the path
(here we identify vertices of $\bt$ and $\bt_{A_n}$).
\begin{Proposition}
\label{P1234}
Let $(\alpha,\beta)$ and $(\gamma,\delta)$ be pairs in  $L(\bt) \times R(\bt)$.
As $n \to \infty$
\[ \Pr_\tau(\bpi(1,2) \mbox{ contains $e_n$ with entrance-exit pair } (\alpha,\beta); \ 
\bpi(3,4) \mbox{ contains $e_n$ with entrance-exit pair } (\gamma,\delta) )
\] \[
\leq (1 + o(1)) \kappa_{\alpha,\beta,\gamma,\delta} 
\left[ \sfrac{\log n}{n} e^{-\sigma} \right]^2
\]
where 
\begin{equation}
\kappa_{\alpha,\beta,\gamma,\delta} = 2^{\ind \{\gamma = \alpha\}
+ \ind \{\delta = \beta\}}. 
\label{kappa-def}
\end{equation}
\end{Proposition}
Proposition \ref{P1234} (more precisely, the variant Proposition \ref{P1234-hack} described in section \ref{sec-hack}) will be proved in section \ref{sec-var}.
Intuitively we have ``$= (1 + o(1))$" instead of ``$\leq (1+o(1))$", but the fact that we need only prove an inequality is technically helpful.  Also intuitively, the constant $\kappa$ arises as 
\[ \kappa_{\alpha,\beta,\gamma,\delta} = \Ex [ W_\alpha W_\beta W_\gamma W_\delta] \]
where $W_\alpha$ is the limit Exponential($1$) r.v.  arising 
(cf. Lemma \ref{Lkt}) in the growth of the percolation process from source $\alpha$ using flows avoiding any other vertex of the neighborhood.

{\bf Proof of Proposition \ref{Pcond-var}.} 
The sum of $\kappa_{\alpha,\beta,\gamma,\delta}$ over all choices of  $\alpha,\beta,\gamma,\delta$ works out as $\# L(\bt) \# R(\bt) (\# L(\bt) +1)( \#R(\bt) +1)$.
So %(all probabilities conditional on $\NN_\tau(n-1,n) = \bt_{A_n}$)
 \[ \Pr_\tau(\bpi(1,2) \mbox{ contains $e_n$ }; \ 
\bpi(3,4) \mbox{ contains $e_n$ } )
\leq (1 + o(1)) \# L(\bt) \# R(\bt) (\# L(\bt) +1)( \#R(\bt) +1)
\left[ \sfrac{\log n}{n} e^{-\sigma} \right]^2 .
\]
Using Corollary \ref{epi} we get a covariance bound 
\[  \Pr_\tau(\bpi(1,2) \mbox{ contains $e_n$ }; \ 
\bpi(3,4) \mbox{ contains $e_n$ } ) 
- \Pr_\tau^2(\bpi(1,2) \mbox{ contains $e_n$ })
\]
 \begin{equation}
\leq (1 + o(1)) \# L(\bt) \# R(\bt) (\# L(\bt) +  \#R(\bt) +1)
\left[ \sfrac{\log n}{n} e^{-\sigma} \right]^2 . 
\label{covar1234}
\end{equation}
The contribution to $F_n(e_n)$ from source-destination pairs $(i,j)$ where $i$ or $j$ is in $\bt$ is negligible, so we may replace $F_n(e_n)$ by 
\[ G_n(e_n) := \sfrac{1}{n} \sum_{(i,j)} \ind \{e_n \in
\bpi(i,j)\} \]
where here and below the sum is over ordered pairs of vertices in $[n] \setminus A_n$.  
Writing $e = e_n$,  
\[ \var_{\!\! \tau} \sum_{(i,j)} \ind \{e \in \bpi(i,j)\} 
\leq \Ex_\tau  \sum_{(i,j)} \ind \{e \in \bpi(i,j)\} \]
\[
+ \sum_{(i,j)} \sum_{(i^\prime,j^\prime) \neq (i,j)} 
\left[  \Pr_\tau(\bpi(i,j) \mbox{ contains $e$ }; \ 
\bpi(i^\prime,j^\prime) \mbox{ contains $e$ } ) 
- \Pr_\tau^2(\bpi(i,j) \mbox{ contains $e$ })
\right] . \]  
Using symmetry and a compatibility condition 
(a directed edge $e$ cannot be in the shortest path from $i$ to $j$ and also in the shortest path from $j$ to $k$) we find
\[ \var_{\!\! \tau} G_n(e) \leq n^{-1} \Ex_\tau G_n(e) 
+ 2n \Pr_\tau(\bpi(1,2) \mbox{ contains $e$ }; \ 
\bpi(1,3) \mbox{ contains $e$ } ) \]
\[ + n^2 
\left[
\Pr_\tau(\bpi(1,2) \mbox{ contains $e$ }; \ 
\bpi(3,4) \mbox{ contains $e$ } ) 
- \Pr_\tau^2(\bpi(1,2) \mbox{ contains $e$ }) 
\right] . \]   
The first term is $O(n^{-1} \log n)$ by Proposition \ref{Pcond-mean}.  
Bounding the second term crudely by 
$2n \Pr_\tau(\bpi(1,2) \mbox{ contains $e$ })$
 and using Corollary \ref{epi} shows the second term is $O( \log n)$.  
 So the dominant term is the third term, which by (\ref{covar1234}) shows 
 \[ \var_{\!\! \tau} G_n(e) 
\leq (1 + o(1)) \# L(\bt) \# R(\bt) (\# L(\bt) +  \#R(\bt) +1)
\left[ \log n \  e^{-\sigma} \right]^2 . \]
Combining with Proposition \ref{Pcond-mean} we have established Proposition \ref{Pcond-var}.
\qed

\subsection{WLLN for a local functional}
\label{sec-WLLN}
The point of Propositions \ref{Pcond-mean} and \ref{Pcond-var} 
is that the normalized flow $F_n(e)/\log n$, which {\em a priori} involves the global structure of $\GG_n$,  can be approximated by a certain functional,
$\phi_n^\tau(e)$ below, which depends only on the ``local" structure of $\GG_n$ near $e$.
It is a general fact that empirical (random) distributions of such ``local functionals" on 
$\GG_n$ converge to the limit non-random distribution associated with the  Yule processs/PWIT mentioned in section \ref{sec-Yule}.
Rather than prove a general result in this context 
(for the general result in a different context see 
\cite{me52} Proposition 7)
we will just derive the specific result we need,
Proposition \ref{Pphitau}.

Fix $\tau > 0$. 
Recall from section \ref{sec-Ls} the definition of the neighborhood 
$\NN_\tau (e)$ of a directed edge $e$ of $\GG_n$.
For each directed edge $e$ define 
\begin{equation} \phi^\tau_n(e) = \# L(\NN_\tau(e))  \# R(\NN_\tau(e)) \exp(-2 \tau - L_e) 
\ind \{L_e \leq \tau\} \ind\{\ \NN_\tau(e) \mbox{ is a tree}\}. \label{phitau-def} 
\end{equation}
Define $\Phi^\tau_n$ as the empirical measure on $[0,\infty)^2$ obtained by putting weight $1/n$  on each point 
$(L_e,\phi^\tau_n(e))$ associated with the edges $e$ of $\GG_n$ for which $L_e \leq \tau$:
\[ \Phi^\tau_n(\cdot,\cdot) = \sfrac{1}{n} \sum_{\mbox{\tiny{directed }} e} \ind \{ (L_e, \phi^\tau_n(e)) \in (\cdot,\cdot)\}  \]
and define the mean measure 
\[ \bar{\Phi}_n^\tau(\cdot,\cdot) = \Ex \Phi^\tau_n(\cdot,\cdot) . \]
Define a limit measure 
\[
\bar{\Phi}_\infty^\tau(\cdot,\cdot) = 
\int_0^\tau \Pr \left((u,\sfrac{W^\tau_1 W^\tau_2}{e^{2\tau}} e^{-u}) \in \cdot \times \cdot \right) \ du \]
where $W^\tau_1$ and $W^\tau_2$ are independent Geometric($e^{-\tau}$).
So $\bar{\Phi}^\tau_\infty$ has total mass $\tau$.
\begin{Proposition}
\label{Pphitau}
For any continuous test function 
$h:[0,\infty)^2 \to \Reals$ 
with compact support,
\[ \int h \ d\Phi^\tau_n 
\to \int h \ d\bar{\Phi}^\tau_\infty 
\mbox{ in } L^1 . \]
\end{Proposition}

The proof rests upon the following straightforward lemma.  
Although superficially similar to Proposition \ref{P1234} in using vertices 
$1,2,3,4$ as typical vertices, their role here is different.
The precise statement is a bit fussy because the neighborhood must be a tree in 
order for left and right sides to be well-defined.
\begin{Lemma}
\label{LtN}
Fix $\tau > 0$ and the directed edges $(1,2)$ and $(3,4)$. 
\\(a) $n \Pr(L_{12} \leq \tau, L_{12} \in \cdot)$ converges vaguely to Lebesgue measure on $[0,\tau]$. 
\\(b) Uniformly on $\{L_{12} \leq \tau\}$ 
\begin{equation}
\Pr(\NN_\tau(1,2) \mbox{ is a tree}|L_{12}) \rightarrow 1
\end{equation}
as $ n\rightarrow \infty $.
The same holds for edges $(1,2)$ and $(3,4)$; that is, uniformly on the set $\{L_{12} \leq \tau\} \cap \{L_{34} \leq \tau \}$
\begin{equation}
\Pr(\{\NN_\tau(1,2) \mbox{ is a tree}\} \cap \{\NN_\tau(3,4) \mbox{ is a tree}\}|L_{12},L_{34}) \rightarrow 1 .
\end{equation} 
\\(c) 
Let $0 < \ell_{12}, \ell_{34} \leq \tau$. 
Write $(\widetilde{N}^\tau_n(1),\widetilde{N}^\tau_n(2))$ 
for the numbers of vertices in the left and right sides of the neighborhood 
$\NN_\tau(1,2)$;
define $(\widetilde{N}^\tau_n(3),\widetilde{N}^\tau_n(4))$ 
similarly for the neighborhood 
$\NN_\tau(3,4)$ 
(these are well-defined when the neighborhoods are trees).
Conditional on the event 
$\{\NN_\tau(1,2) \mbox{ is a tree}\} \cap \{\NN_\tau(3,4) \mbox{ is a tree}\}
\cap \{L_{12} = \ell_{12}, L_{34} = \ell_{34}\}$
we have
\[ (\widetilde{N}^\tau_n(1),\widetilde{N}^\tau_n(2),\widetilde{N}^\tau_n(3),\widetilde{N}^\tau_n(4)) 
\cd
(W^\tau_1,W^\tau_2,W^\tau_3,W^\tau_4)
\]
where the $W$'s are independent Geometric($e^{-\tau}$).
\end{Lemma}

{\bf Proof of Proposition \ref{Pphitau}.}
The mean measure 
$\bar{\Phi}^\tau_n$
equals
\[ \sfrac{1}{n} \times 
n(n-1) \Pr \left(L_{12} \leq \tau, L_{12} \in \cdot,
\frac{\widetilde{N}^\tau_n(1) \widetilde{N}^\tau_n(2)}{e^{2\tau}} e^{-L_{12}} \in \cdot \right)
.\]
Because
$n \Pr (L_{12} \leq \tau, L_{12} \in \cdot)$
converges vaguely to 
Lebesgue measure on $[0, \tau]$,
Lemma \ref{LtN} 
(here only vertices $1$ and $2$ are relevant)
implies vague convergence 
$\bar{\Phi}^\tau_n
 \to 
\bar{\Phi}^\tau_\infty$
of mean measures.
To get $L^2$ convergence 
it is enough to show that for a generic test function $h$ we have 
\[ \var \left( 
\int h \ d 
\Phi^\tau_n
\right) = 
 \var \left( 
\sfrac{1}{n} \sum_{e: L_e \leq \tau} h(L_e, \phi_n(e)) \right) 
\to 0 . \]
Expanding the right side as the variance-covariance sum,
the contribution to variance from terms 
$(e,e^\prime)$ with $4$ distinct end-vertices 
tends to $0$
by Lemma \ref{LtN} and the fact that 
$n^2 \Pr (L_{12} \leq \tau, \NN_\tau(1,2)\mbox{ is a tree}, L_{12} \in \cdot,L_{34} \leq \tau,\NN_\tau(3,4)\mbox{ is a tree }, L_{34} \in \cdot)$
converges vaguely to 
Lebesgue measure on $[0, \tau]^2$.
The contribution from terms with $3$ distinct end-vertices 
is bounded by 
\[ \frac{1}{n} ||h||_\infty^2 
\Ex \Delta_1^2(\tau) \to 0 \]
where $\Delta_1(\tau)$ 
is the number of edges at $1$ with length less than $\tau$.  
And the contribution from pairs $(e,e)$ is bounded by 
$ ||h||_\infty^2 \Pr(L_{12} \leq \tau) \to 0$.

\subsection{Completing the proof of Theorem \ref{T1}}
\label{sec-complete-proof}
The remainder of the proof uses only ``soft" arguments.

Propositions \ref{Pcond-mean} and \ref{Pcond-var} were stated for fixed
 $\bt \in \bT_\tau$, but it is clear that convergence is uniform over subsets of $\bT_\tau$ on which  the length of distinguished edge is bounded and the number of vertices is bounded.  
Rephrasing those Propositions gives, after some obvious manipulations:
\begin{Corollary}
Fix $\tau > 0$ and $K < \infty$.  
As $n \to \infty$
\[
\Ex \left( \left.  \frac{F_n(e_n)}{\log n} \right| \NN_\tau(e_n) \right) 
- \phi^\tau_n(e_n) \to 0 \]
\[ \var \left( \left.  \frac{F_n(e_n)}{\log n} \right| \NN_\tau(e_n) \right) 
\leq 3K^3e^{-\tau} + o(1) \]
uniformly over $e_n$ satisfying 
\[\NN_\tau(e_n) \mbox{ is a tree };\quad L_{e_n} \leq \tau; \quad 
\max( \# L(\NN_\tau(e_n)) , \# R(\NN_\tau(e_n))) \leq K e^\tau . \]
\end{Corollary} 
Now fix $\eps > 0$.  
Applying Chebyshev's inequality (conditional on $\NN_\tau(e_n)$) and
taking limits,
\[
\limsup_n 
\Ex 
\sfrac{1}{n} \sum_{\mbox{\tiny{directed }} e} 
\ind \{L_e \leq \tau,\NN_\tau(e) \mbox{ is a tree}\} 
\ind \{\max( \# L(\NN_\tau(e)) , \# R(\NN_\tau(e))) \leq K e^\tau\} 
\ind \{ | \sfrac{F_n(e)}{\log n} - \phi^\tau_n(e)| >  \eps \} \]
 \[ \leq 3 \eps^{-2} K^3 e^{- \tau} . \]
And using Lemma \ref{LtN}
\begin{eqnarray*}
\lefteqn{ \limsup_n 
\Ex 
\sfrac{1}{n} \sum_{\mbox{\tiny{directed }} e} 
\ind \{L_e \leq \tau\} 
\ind \{\NN_\tau(e) \mbox{ is a tree }\}
\ind \{\max( \# L(\NN_\tau(e)) , \# R(\NN_\tau(e))) > K e^\tau\} 
}\\
&=& 
\lim_n n \Ex \ind \{L_{12} \leq \tau\} 
\ind \{\NN_\tau(1,2) \mbox{is a tree }\}
\ind \{\max( \# L(\NN_\tau(1,2)) , \# R(\NN_\tau(1,2))) > K e^\tau\} \\
&=& \tau \Pr (\max(W^\tau_1,W^\tau_2) > Ke^\tau)\\
&\leq & 2 \tau \exp(- K + e^{- \tau}) .
\end{eqnarray*}
Combining these bounds, 
\[   \limsup_n 
\Ex 
\sfrac{1}{n} \sum_{\mbox{\tiny{directed }} e} 
\ind \{L_e \leq \tau\} \ind \{\NN_\tau(e) \mbox{ is a tree }\} \ind \{ | \sfrac{F_n(e)}{\log n} - \phi^\tau_n(e)| >  \eps \} 
\leq 3 \eps^{-2} K^3 e^{- \tau} + 2 \tau \exp(- K + e^{- \tau}) . \]
Apply this with $K = \tau$ and then let $\tau \to \infty$:
\begin{equation}
  \lim_\tau  \limsup_n 
\Ex 
\sfrac{1}{n} \sum_{\mbox{\tiny{directed }} e} 
\ind \{L_e \leq \tau\} \ind \{\NN_\tau(e) \mbox{ is a tree }\} \ind \{ | \sfrac{F_n(e)}{\log n} - \phi^\tau_n(e)| >  \eps \} 
= 0 . \label{tau-n}
\end{equation}

Also by Lemma \ref{LtN} for each fixed $\tau$

\begin{equation}
\label{tau-not-tree}
\limsup_n \Ex \sfrac{1}{n} \sum_{\mbox{\tiny{directed }} e}  \ind\{ L_e\leq \tau , \NN_\tau(e) \mbox{ not a tree}\} = 0\end{equation}

We now want to be a little fussy about the underlying space for our bivariate measures, which we will take to be 
$[0,\infty) \times (0,\infty)$ 
(recall the first coordinate is length, the second is flow).  This means that the $\sigma$-finite limit measure $\psi$ arising in the statement of Theorem \ref{T1} is finite on compact subsets.
Recall that 
vague convergence $\nu_n \to \nu$ of measures on $[0,\infty) \times (0,\infty)$ 
means
$\int h d\nu_n \to \int h d\nu$ 
for bounded continuous test functions $h:[0,\infty) \times (0,\infty) \to \Reals$ with {\em compact} support, 
and that in checking vague convergence we need consider only test functions with finite
Lipschitz norm $|| h ||_\Lip$.  
A random measure can be viewed as a random variable taking values in the space of measures equipped with the vague topology, and so it makes sense to consider convergence in probability
\begin{equation}
\psi_n \to \psi \mbox{ in probability } 
\label{psi-psi}
\end{equation}
for the random measures $\psi_n$ appearing in Theorem \ref{T1}, and this is what we shall prove.  
Of course it suffices to prove that for test functions $h$ we have convergence in probability for the $\Reals$-valued random variables $ \int h d\psi_n$, and we shall prove the stronger result 
\begin{equation}
\int h d\psi_n \to \int h d\psi \mbox{ in } L^1 . 
\label{hpsih}
\end{equation}
To prove this, recall the definitions 

\[ \psi_n(\cdot,\cdot) = \sfrac{1}{n} \sum_{\mbox{\tiny{directed }} e} \ind \{ (L_e, F_n(e)/\log n) \in (\cdot,\cdot)\}  \]
\[ \Phi^\tau_n(\cdot,\cdot) = \sfrac{1}{n} \sum_{\mbox{\tiny{directed }} e} \ind \{ (L_e, \phi^\tau_n(e)) \in (\cdot,\cdot)\}   . \]
Fix $h$ with support contained in $[0,\tau_0] \times [0,\infty)$.
Then for $\tau > \tau_0$
\begin{eqnarray*}
\left| \int h d \psi_n - \int h d \Phi^\tau_n \right | &\leq& 
2 || h ||_\infty 
\sfrac{1}{n} \sum_{\mbox{\tiny{directed }} e} 
\ind \{L_e \leq \tau_0, \NN_\tau(e) \mbox{ is a tree }\} \ind \{ | \sfrac{F_n(e)}{\log n} - \phi^\tau_n(e)| >  \eps \}  \\ 
&+& 2 || h ||_\infty
\sfrac{1}{n} \sum_{\mbox{\tiny{directed }} e}
\ind\{L_e \leq \tau_0, \NN_\tau(e) \mbox{ not a tree }\} +
  \eps || h ||_\Lip 
\sfrac{1}{n} \sum_{\mbox{\tiny{directed }} e} 
\ind \{L_e \leq \tau_0\} . 
\end{eqnarray*}
Because 
$\Ex \sfrac{1}{n} \sum_{\mbox{\tiny{directed }} e} 
\ind \{L_e \leq \tau_0\} \to \tau_0$ 
we can use (\ref{tau-n},\ref{tau-not-tree}) to deduce
\[  \lim_\tau  \limsup_n 
\Ex \left| \int h d \psi_n - \int h d \Phi^\tau_n \right | \leq  \eps || h ||_\Lip \tau_0 . \] 
Because $\eps$ is arbitrary, this shows 
\[  \lim_\tau  \limsup_n 
\Ex \left| \int h d \psi_n - \int h d \Phi^\tau_n \right | = 0. \]  
Proposition \ref{Pphitau} allows us to replace the random measure $\Phi^\tau_n $ 
by the limit mean measure $\bar{\Phi}^\tau_\infty$: 
 \[  \lim_\tau  \limsup_n 
\Ex \left| \int h d \psi_n - \int h d \bar{\Phi}^\tau_\infty  \right | = 0. \]   
But $\bar{\Phi}^\tau_\infty \to \psi$ vaguely as $\tau \to \infty$, 
and so we have proved (\ref{hpsih}) and thence (\ref{psi-psi}),
which is our formalization of the final assertion of Theorem \ref{T1}.

To prove the other assertion (\ref{T1-Ex}) of Theorem \ref{T1}, 
recall that the fact (\ref{Dn}) 
$\Ex \bar{D}_n \sim \log n$
becomes, via (\ref{ell-y}), 
\[ 
\Ex \int \ell y \  \psi_n(d \ell, dy) \to 
\int \ell y \ \psi(d\ell, dy) = 1 . \]  
This enables us to extend the $L^1$ convergence (\ref{hpsih}) 
from continuous $h$ with compact support 
to continuous $h \geq 0$ satisfying 
$\sup_{\ell, y} \frac{h(\ell,y)}{\ell y} < \infty$.  
Using such functions to approximate
the function 
$\ind \{\ell > \eps, y > z\}$ 
shows
\[ \sfrac{1}{n} \# \{e: L_e > \eps, F_n(e) > z \log n\}
\to_{L^1}
\psi( (\eps,\infty) \times (z,\infty)) . \] 
Because 
$\Ex \sfrac{1}{n} \# \{e: L_e \leq  \eps\} \to \eps$ 
we can let $\eps \to 0$ and deduce
\[ \sfrac{1}{n} \# \{e:  F_n(e) > z \log n\}
\to_{L^1}
\psi( (0,\infty) \times (z,\infty))  \] 
which is the first assertion of Theorem \ref{T1}.

\subsection{Distance based truncation of flows}
\label{sec-hack}
To avoid notational complications, the exposition above omitted one technical point.  
Recall that path-lengths $D(i,j)$ are $\log n \pm O(1)$ in probability.
In seeking to prove Proposition \ref{P1234} there are technical difficulties with unusually long paths, 
which we will handle by truncating them out.  
Precisely, instead of proving Proposition \ref{P1234} we will prove 
\begin{Proposition}
\label{P1234-hack} 
Fix $B < \infty$.
Let $(\alpha,\beta)$ and $(\gamma,\delta)$ be pairs in  $L(\bt) \times R(\bt)$.
Conditionally on $\NN_\tau(n-1,n) = \bt_{A_n}$, as $n \to \infty$
\[ \Pr(\bpi(1,2) \mbox{ contains $e_n$ with entrance-exit pair } (\alpha,\beta);  
\ \len(\bpi(1,2)) \leq \log n + B; \] \[
\bpi(3,4) \mbox{ contains $e_n$ with entrance-exit pair } (\gamma,\delta) , \ 
\ \len(\bpi(3,4)) \leq \log n + B)
\] \[
\leq (1 + o(1)) \kappa_{\alpha,\beta,\gamma,\delta} 
\left[ \sfrac{\log n}{n} e^{-\sigma} \right]^2
\]
where 
\begin{equation}
\kappa_{\alpha,\beta,\gamma,\delta} = 2^{\ind \{\gamma = \alpha\}
+ \ind \{\delta = \beta\}}. 
\label{kappa-def-2}
\end{equation}
\end{Proposition}
In this section we explain 
(omitting some details at the end) 
why it is enough to prove Proposition \ref{P1234-hack} in place of Proposition \ref{P1234}.
Consider the analog of flow $F_n(e)$ when contributions from paths of length 
$> \log n + B$ are ignored:
\[ F^{[B]}_n(e) := \sfrac{1}{n} \sum_{i \in [n]} \sum_{j \in [n], j \neq i} \ind \{e \in
\bpi(i,j)\} \ind \{\len (\bpi(i,j)) \leq \log n + B\}  \leq F_n(e) .\]
An easy argument shows that for large $B$ the global effect of truncation is negligible:
\begin{Lemma}
\label{Ltruncate}
\[ \lim_{B \to \infty}  \limsup_n 
\Ex \sfrac{1}{n} \sum_e L_e \sfrac{F_n(e) - F^{[B]}_n(e)}{\log n} 
= 0 . \]
\end{Lemma}
\proof 
Recall $D(i,j) = \len(\bpi(i,j))$.
Calculating the effect of truncation on edge-flows and on source-destination distances gives the identity 
\[ \sum_e L_e (F_n(e) - F^{[B]}_n(e)) = \sfrac{1}{n} \sum_{j \neq i} D(i,j) \ind \{D(i,j) > \log n + B\} . \]  
Taking expectations and using symmetry 
\[ \Ex \sfrac{1}{n} \sum_e L_e (F_n(e) - F^{[B]}_n(e)) 
= \sfrac{n(n-1)}{n^2} \Ex D(1,2) \ind \{D(1,2) > \log n + B\} . \]  
The result now follows from the mean and variance limits at (\ref{D12-limits}).
\qed

Now we can choose $B_n \uparrow \infty$ sufficiently slowly that 
(by Lemma \ref{Ltruncate})
\begin{equation}
\Ex \sfrac{1}{n} \sum_e L_e \sfrac{F_n(e) - F^{[B_n]}_n(e)}{\log n} 
\to 0 \label{L-trunc-2}
\end{equation}
and such that 
\[
\mbox{the conclusion of Proposition \ref{P1234-hack} holds for $B_n$ in place of $B$} . 
\]
The idea is now to repeat the arguments in sections
\ref{sec-cv} - \ref{sec-complete-proof} 
using the truncated flow 
$F_n^{[B_n]}(e_n)$
in place of $F_n(e_n)$.
This will establish Theorem \ref{T1} for the truncated flows, but then
(\ref{L-trunc-2}) establishes it for untruncated flows.
The arguments would go through unchanged if we knew the truncated version
of the conditional mean estimate of Proposition \ref{Pcond-mean}:
\begin{equation}
\expec_\tau(F_n^{[B_n]}(e_n)) = (1+o(1)) \#L(\bt) \#R(\bt) e^{-\sigma} \log{n} . 
\label{P11-trunc}
\end{equation}
Of course the 
conditional 
upper bound
``$\leq (1+o(1))$" in (\ref{P11-trunc})
follows from Proposition \ref{Pcond-mean}, but we need the lower bound 
``$\geq (1+o(1))$" in (\ref{P11-trunc})
in order to go from the upper bound on second moment to the upper bound on variance 
 -- cf. (\ref{covar1234}).
However, from the conditional upper bound and because 
(\ref{L-trunc-2})
implies a lower bound for unconditional expectation,
Markov's inequality
implies that the conditional lower bound 
in (\ref{P11-trunc}) holds
for all neighborhoods 
$\NN_\tau(n-1,n)$ 
excluding some occuring with probability $\to 0$ as $n \to \infty$.
And this is enough to complete the proof of Theorem \ref{T1}.

\section{The variance estimate}
\label{sec-var}
This section is devoted to the proof of Proposition \ref{P1234-hack},
which will complete the proof of Theorem \ref{T1}.
Let us repeat the ``conditioned" setting that we work in,
throughout the section.
There is a fixed tree $\bt$ 
with distinguished directed edge $e$.
In the network $\GG_n$ we fix edge 
$e_n = (n-1,n)$
and label set 
$A_n = \{n - \#\bt +1, \ldots, n\}$.
Label $\bt$ as $\bt_{A_n}$
so that $e$ is labeled $(n-1,n)$. 
Then condition on 
$\NN_\tau(n-1,n) = \bt_{A_n}$. 
Recall that Lemma \ref{Lneigh} tells us the effect of this conditioning.
 In particular, for $i \not\in A_n, \ j \in A_n$ the length of the edge-segment $(i,j)$ from $i$ to the boundary of the neighborhood has Exponential($1/n$) distribution, independently for different edges.

Roughly speaking, the issue in proving Proposition \ref{P1234-hack} is to estimate the dependence
between the events
\\ (i)
the shortest path $\bpi(1,2)$ between vertex $1$ and vertex $2$ uses edge $e_n$
\\
(ii) 
the shortest path $\bpi(3,4)$ between vertex $3$ and vertex $4$ uses edge $e_n$.

\noindent
Corollary \ref{epi} tells us the asymptotic probabilities 
of these events, so a natural approach is to condition on (i) 
and seek to calculate the conditional probability of (ii).
Now (i) breaks into two assertions:
\\
(ia) there is a short path 
(length $s = \log n \pm O(1)$)
from $1$ to $2$ via $e$;
\\
(ib) there is no shorter path from $1$ to $2$.

Now conditioning on (ia) can be implemented by conditioning on all edges in the path, and this doesn't affect lengths of other edges of $\GG_n$.
But event (ib) implicitly specifies that alternate short routes do not exist, and the effect of this conditioning on 
other edge-lengths of $\GG_n$ (while intuitively small) 
seems hard to handle rigorously. 
Instead, we shall carefully organize an argument to avoid ever conditioning on any 
``shortest path" event.
In outline, the argument has three steps.
\begin{itemize}
\item Calculate chance of existence of paths $\pi_{12}, \pi_{34}$ of specified lengths through $e_n$
(section \ref{joint-intensity})
\item Conditional on existence of such paths, what is the chance they are the {\em shortest} paths?  
The percolation processes from vertices $i = 1,2,3,4$ avoiding edge $e_n$  become approximately 
{\em size-biased} Yule processes (section \ref{sec-sizebias}) reaching $\widetilde{W}_i e^t$ vertices in distance $t$; so the chance of a path from $1$ to $2$ of length $t$ avoiding $e_n$ is approximately 
$\exp(- \widetilde{W}_1 \widetilde{W}_2 e^t)$ (section \ref{sec-P22}). 
\item These two estimates are combined in section \ref{cond-distr-other-short}.
\end{itemize}

\subsection{Joint intensity for two short paths through $e$}
\label{joint-intensity}
For this section we introduce some handy notation.  
We will describe the relationship 
(for an event $B$, a random variable $T$, and a function $f$)
\[ \Pr(B, \ T \in (t,t+dt)) = f(t) dt \]
by the phrase
\[ \mbox{the event $[[B, \ T=t ]]$ has intensity $f(t)$.} \]
But we will describe events in words, rather than inventing ad hoc symbols, using the  
brackets $[[ \ldots \ldots ]]$ to highlight the verbal description of the event.

For $\alpha \in L(\bt)$ and $\beta \in R(\bt)$ 
define 
$h_{\alpha,\beta}(s_1,s_2) $
to be the intensity (in $s_1$ and $s_2$) of the event:
\begin{quote}
[[there exists a path from $1$ to $v_L$ which crosses the neighborhood boundary 
at time $s_1$ and then first hits vertex $\alpha$;\\
and
there exists a path from $2$ to $v_R$ which crosses the neighborhood boundary 
at time $s_2$ and then first hits vertex $\beta$.]]
\end{quote}
Note that such paths, linked via the path from $\alpha$ to $\beta$ in the neighborhood,
specify a path $\pi_{12}$ of length $s_1 + s_2 + \sigma$ from $1$ to $2$ via $e$.
This path may or may not be the shortest path from $1$ to $2$, 
depending on lengths of other edges in $\GG_n$.

Given also $\gamma \in L(\bt)$ and $\delta \in R(\bt)$, define an event which replicas the event above:
\begin{quote}
[[there exists a path from $3$ to $v_L$ which crosses the neighborhood boundary 
at time $t_1$ and then first hits vertex $\gamma$;\\
and
there exists a path from $4$ to $v_R$ which crosses the neighborhood boundary 
at time $t_2$ and then first hits vertex $\delta$.]]
\end{quote}
Again, such paths specify a path $\pi_{34}$ from $3$ to $4$ via $e$.
In order for it to be possible that \underline{both} $\pi_{12}$ and $\pi_{34}$ are shortest paths, the following simple {\em compatibility conditions} must hold.\\
(i) if the paths from $1$ and from $3$ meet at some vertex $v_*$ outside the neighborhood, 
then they must coincide from $v_*$ to the neighborhood. \\
(ii) if the paths from $2$ and from $4$ meet at some vertex $v^*$ outside the neighborhood, 
then they must coincide from $v^*$ to the neighborhood. \\ 
(iii) the set of vertices visited by the paths from $1$ and $3$ must be disjoint from 
the set of vertices visited by the paths from $2$ and $4$.

Define 
$H_{\alpha,\beta,\gamma,\delta}(s_1,s_2,t_1,t_2) $
to be the intensity of 
both events happening and the compatibility conditions holding.
\begin{Lemma}
\label{Lkappa}
\begin{eqnarray}
h_{\alpha,\beta}(s_1,s_2) 
&\leq& \frac{1}{n^2} \exp(s_1 + s_2 ) 
\label{hab}\\
H_{\alpha,\beta,\gamma,\delta}(s_1,s_2,t_1,t_2)  
&\leq& \frac{\kappa }{n^4} \exp(s_1 + s_2 +t_1 + t_2 ) .
\label{Habcd}
\end{eqnarray}
Here 
$\kappa =  \kappa_{\alpha,\beta,\gamma,\delta}$ as at (\ref{kappa-def}). 
\end{Lemma}
Note these are inequalities for finite $n$. 
Heuristically they are asymptotic equalities in the ranges of interest to us.

 \proof
The argument is based on exact formulas, starting with the following.  
The intensity of the event 
\begin{quote}
 [[there exists a path from $1$ to $v_L$ which crosses the neighborhood boundary 
at time $s_1$ and then first hits vertex $\alpha$, taking exactly $k+1$ steps]]
\end{quote}
\begin{equation}
 = \frac{(n-1- \# \bt)_k}{n^k}\  \frac{s_1^k}{k!} \  \ \sfrac{1}{n} \ \exp( \sfrac{-s_1 - b(\alpha)}{n}) 
 \label{1vL}
 \end{equation}
where we recall that $b(\alpha)$ is the distance from $\alpha$ to the neighborhood boundary.  
To prove this, take 
$0<u_1 < u_2 < \ldots < u_k < s_1$ 
and consider the probability that the $j$'th step ends at distance $[u_j, u_j + du_j]$ from $1$ 
and the boundary crossing is at distance $[s_1,s_1 + ds_1]$.  
This probability is
\begin{equation} (n - 1 - \# \bt)_k \times \prod_{j=1}^k \sfrac{1}{n} \exp(-(u_j - u_{j-1})/n) du_j 
\times \sfrac{1}{n} \exp(-(s_1 + b(\alpha) - u_k)/n) ds_1 \label{intens-1} 
\end{equation}
where the first term indicates number of choices of $k$ intermediate vertices, and the other terms 
are the Exponential(mean $n$) densities of edge-lengths.  
Because 
\begin{equation} \mbox{$\int \ldots \int $}_{0<u_1 < u_2 < \ldots < u_k < s_1} \ du_1 \ldots du_k 
= s_1^k/k! 
\label{kfold}
\end{equation}
we deduce (\ref{1vL}).  

The first and last terms of (\ref{1vL}) are $ \leq 1$.
Summing over $k$ shows that the intensity of 
\begin{quote}
[[there exists a path from $1$ to $v_L$ which crosses the neighborhood boundary 
at time $s_1$ and then first hits vertex $\alpha$]]
\end{quote}
is $\leq \frac{1}{n} \exp(s_1)$.  
Combining this with the similar argument on the right side of the neighborhood 
gives (\ref{hab}).  

To prove (\ref{Habcd}), 
because the left and right sides have analogous arguments,
the issue is to show that the intensity of the event
\begin{quote}
[[there exists a path from $1$ to $v_L$ which crosses the neighborhood boundary 
at time $s_1$ and then first hits vertex $\alpha$;\\
there exists a path from $3$ to $v_L$ which crosses the neighborhood boundary 
at time $t_1$ and then first hits vertex $\gamma$]]
\end{quote}
\begin{equation}
\leq 2^{\ind (\alpha = \gamma)} 
\frac{1}{n^2} \exp(s_1 + t_1). 
\label{1v3}
\end{equation}

Now in the case $\alpha \neq \gamma$, or for the contribution to the case 
$\alpha = \gamma$ from disjoint paths, we get  density 
$\leq \frac{1}{n^2} \exp(s_1 + t_1)$
by arguments analogous to above.  
Let us show details of the more interesting case 
where $\alpha = \gamma$ and we consider the contribution from merging paths.
Consider the intensity (in $r, s_1, t_1$) of the event:
\begin{quote}
[[there exists a path from $1$ to $v_L$ which crosses the neighborhood boundary 
at time $s_1$ and then first hits vertex $\alpha$, taking exactly $k_1 + k_2 +1$ steps;\\
there exists a path from $3$ to $v_L$ which crosses the neighborhood boundary 
at time $t_1$ and then first hits vertex $\alpha$, taking exactly $k_3 + k_2 +1$ steps;\\ 
these paths merge at some vertex $v_*$ at distance $r$ from the neighborhood boundary, 
the path from $v_*$ to $\alpha$ using $k_2 + 1$ steps.]]
\end{quote}
Analogous to (\ref{1vL}) this intensity has an exact formula
 \[ \frac{(n-1- \# \bt)_{k_1+k_2+k_3 -1}}{n^{k_1+k_2+k_3}}\  \frac{(s_1-r)^{k_1}}{ k_1!} \frac{(t_1-r)^{k_3}}{k_3!} \frac{r^{k_2}}{k_2!}  \  \ \sfrac{1}{n} \ \exp( \sfrac{-s_1 -t_1 + r - b(\alpha)}{n}) 
 . \]
This intensity is bounded by
\[ \frac{1}{n^2} 
\frac{(s_1-r)^{k_1} (t_1-r)^{k_3} r^{k_2}}{k_1! k_3! k_2!} 
. \]
Summing over $(k_1,k_2,k_3)$ shows that the intensity of
\begin{quote}
[[there exist paths from $1$ (resp. $3$) to $v_L$ which cross the neighborhood boundary 
at time $s_1$ (resp. $t_1 $) and then first hit vertex $\alpha$, 
having merged at distance $r$ before the boundary]]
\end{quote}
is bounded by
$\frac{1}{n^2} \exp(s_1 +t_1 -r)$.
Integrating over $r$ shows that the contribution to (\ref{1v3}) 
from merging paths
is $\leq 
\frac{1}{n^2} \exp(s_1 + t_1)$.
This establishes (\ref{1v3}).

\subsection{A Cox point process}
\label{sec-Cox}
Here we introduce a process and a lemma; how the process arises will be seen in section \ref{sec-P22}.

Take independent random variables $\wW_1, \wW_2$ 
with probability density $w e^{-w}$ on $0<w<\infty$.  
Consider the Cox point process defined by: conditional on $\wW_1, \wW_2$ the points form a Poisson process of rate $(\wW_1 \wW_2 e^s, \ - \infty < s < \infty)$.  
%The next lemma gives two easy calculations with this Cox process.  
Let $\wL$ be the position of the leftmost point of this Cox process. 
%Let $I_k$ be the $k$-point intensity function:
%\[ \Pr (\mbox{ some point in $[s_i,s_i+ds_i], \ 1 \leq i \leq k$}) 
%= I_k(s_1,\ldots,s_k) \ ds_1 \ldots ds_k . \]
\begin{Lemma}
\label{LCox}
% $I_k(s_1,\ldots,s_k) = ( (k+1)!)^2 \ \exp( \sum_i s_i)$.\\
$\int_{-\infty}^\infty e^s \Pr (\wL > s) \ ds = 1 $.
\end{Lemma}
\proof 
%For (a) note 
%$ I_k(s_1,\ldots,s_k) = \Ex \prod_{i=1}^k (\wW_1 \wW_2 e^{s_i})$
%and
%$\Ex \wW^k = (k+1)!$.
Note 
$\Pr (\wL > s|\wW_1, \wW_2) = \exp(-\wW_1 \wW_2 e^s)$
and so
\begin{eqnarray*}
\int_{-\infty}^\infty e^s \Pr (\wL > s) \ ds 
&=& \int_{-\infty}^\infty e^s \Ex \exp(-\wW_1 \wW_2 e^s) \ ds\\
&=& \int_0^\infty \Ex \exp(-\wW_1 \wW_2 u) \ du \quad \mbox{ setting } u = e^s \\
&=& \Ex \left[ \frac{1}{\wW_1 \wW_2} \right] 
= \left[ \Ex \frac{1}{\wW_1} \right]^2 = 1 .
\end{eqnarray*}

\subsection{Conditional distributions of other short routes}
\label{cond-distr-other-short}
Recall we are conditioning on $ \NN_\tau(n-1,n) = \bt_{A_n}$, though this is not indicated in notation.

Fix times $s_1, s_2, t_1, t_2$
and vertices (maybe the same) 
$\alpha, \gamma \in L(\bt_{A_n})$
and 
$\beta, \delta \in R(\bt_{A_n})$.
Recall from section \ref{joint-intensity} that 
$H_{\alpha,\beta,\gamma,\delta}(s_1,s_2,t_1,t_2)$
denotes the intensity of the event
\begin{quote}
[[there exists a path $\pi_{12}$
from $1$ to $2$ via $e_n$, where the path from $1$
crosses the boundary 
of the neighborhood at time $s_1$ and then first hits vertex $\alpha$, 
while the reverse path from $2$
crosses the boundary 
of the neighborhood at time $s_2$ and then first hits vertex $\beta$;  
similarly there exists a path $\pi_{34}$
from $3$ to $4$ via $e_n$, where the path from $3$
crosses the boundary 
of the neighborhood at time $t_1$ and then first hits vertex $\gamma$, 
while the reverse path from $4$
crosses the boundary 
of the neighborhood at time $t_2$ and then first hits vertex $\delta$  
]]
\end{quote}
together with certain compatability conditions.
We will write 
$(\ \cdot \ |s_1,s_2,t_1,t_2)$ 
to denote conditioning on this event.
For such paths we have 
\[ \len(\pi_{12}) = s_1 + s_2 + \sigma, \ 
\ \len(\pi_{34}) = t_1 + t_2 + \sigma, \]
and we write $\len(\pi_{12})$ and $\len(\pi_{34})$ for these sums.

Write $S_{12}$ (resp. $S_{34}$) for the length of the shortest path from $1$ to $2$ (resp. from $3$ to $4$) that does not use edge $e$.
Let us first show that Proposition \ref{P1234-hack}
reduces  to the following proposition.
\begin{Proposition}
\label{S1234}
Fix $B < \infty$. 
Uniformly on 
$\{\max( \len(\pi_{12}), \len(\pi_{34})) \leq \log n + B\}$, 
as $n \to \infty$
\[ \Pr (S_{12} > \len (\pi_{12}), \ S_{34} > \len (\pi_{34})| \  s_1,s_2,t_1,t_2)
\leq \Pr (\xi^{12}_1 > \len (\pi_{12}) - \log n, \ \xi^{34}_1 > \len (\pi_{34}) - \log n) 
+ o(1) \]
where 
$(\xi^{12}_1,\xi^{12}_2,\ldots)$ and $(\xi^{34}_1,\xi^{34}_2,\ldots)$ are independent Cox point processes
as described in section \ref{sec-Cox}.
\end{Proposition}
We will prove this in sections \ref{sec-sizebias} - \ref{sec-P22}, but let us first
show how to deduce Proposition \ref{P1234-hack} from Lemma \ref{Lkappa} and 
Proposition \ref{S1234}.

\noindent
{\bf Proof of Proposition \ref{P1234-hack}.}
A path $\pi_{12}$ via $e_n$ using entrance-exit pair $(\alpha,\beta)$ with $\len(\pi_{12}) = s_1 + \sigma + s_2$ is created  as in the definition of  $h_{\alpha,\beta}(s_1,s_2)$ from two paths with lengths-to-boundary  $s_1$ and $s_2$.
Create $\pi_{34}$ similarly, using paths of lengths $t_1$ and $t_2$.  
The quantity in Proposition \ref{P1234-hack} 
\[ \Pr(\bpi(1,2) \mbox{ contains $e_n$ with entrance-exit pair } (\alpha,\beta);  
\ \len(\bpi(1,2)) \leq \log n + B; \] \[
\bpi(3,4) \mbox{ contains $e_n$ with entrance-exit pair } (\gamma,\delta) , \ 
\ \len(\bpi(3,4)) \leq \log n + B)
\]
can be calculated in terms of the intensity $H$ at (\ref{Habcd}) as 
\begin{eqnarray} 
\int \int \int \int 
\ind\{s_1+s_2+\sigma \leq \log n + B\} 
\ind \{t_1+t_2+\sigma \leq \log n + B\} &&  \label{intH} \\
 \Pr(\bpi(1,2) = \pi_{12}, \bpi(3,4) = \pi_{34} |s_1,s_2,t_1,t_2) && \! \! \! \!
H_{\alpha,\beta,\gamma,\delta}(s_1,s_2,t_1,t_2)
\quad ds_1 ds_2 dt_1 dt_2 . \nonumber
\end{eqnarray}
 Now 
 \[ \Pr(\bpi(1,2) = \pi_{12}, \bpi(3,4) = \pi_{34} |s_1,s_2,t_1,t_2) 
 \leq
 \Pr (S_{12} > \len (\pi_{12}), \ S_{34} > \len (\pi_{34})| \  s_1,s_2,t_1,t_2) , \]
 this being an inequality because there might be shorter paths using $e$.
Bounding the right side using Proposition \ref{S1234} gives 
\[ \Pr(\bpi(1,2) = \pi_{12}, \bpi(3,4) = \pi_{34} |s_1,s_2,t_1,t_2) 
\leq  \Pr(\xi^{12}_1 > s^1) \ \Pr( \xi^{34}_1 > t^1) + \eps_{n,B} 
 \]
 where $\lim_n \eps_{n,B} = 0$
and where 
\[ s^1 = s_1+ s_2+\sigma - \log n; \quad 
 t^1 = t_1+ t_2+\sigma - \log n. \] 
To upper bound (\ref{intH}),
first fix $s^1$ and $t^1$ and calculate
\begin{eqnarray*}
&  &
\int_{s_1+s_2+\sigma = \log n + s^1}
\int_{t_1+t_2+\sigma = \log n + t^1} 
H_{\alpha,\beta,\gamma,\delta}(s_1,s_2,t_1,t_2)
 ds_1 dt_1 \\
& \leq & \int_{s_1+s_2+\sigma = \log n + s^1}
\int_{t_1+t_2+\sigma = \log n + t^1} 
\frac{\kappa }{n^4} \exp(s_1 + s_2 +t_1 + t_2 )
 ds_1 dt_1 \mbox{ by  Lemma \ref{Lkappa} } \\
 &=&  \sfrac{\kappa}{n^2} 
 \exp(s^1 + t^1 - 2 \sigma) 
 \times 
 (\log n + s^1 - \sigma)
 (\log n + t^1 - \sigma)  \\
 &\leq & \kappa e^{-2\sigma} \sfrac{(\log n + B)^2}{n^2} \exp(s^1 + t^1)  
 \end{eqnarray*}
where in the final line we assume 
$\max(s^1,t^1) \leq B$.  
Thus the quantity (\ref{intH})
is bounded by
 \begin{eqnarray*}
 &&  \kappa  \sfrac{(\log n + B)^2}{n^2} e^{-2\sigma}
 \int_{-\infty}^B  \int_{-\infty}^B 
 \left(\Pr(\xi^{12}_1 > s^1) \Pr(\xi^{34}_1 > t^1) + \eps_{n,B} \right) \ e^{s^1} ds^1 \ e^{t^1} dt^1  \\
 && \leq  \kappa  \sfrac{(\log n + B)^2}{n^2} e^{-2\sigma}
 \ (1 + \eps_{n,B} e^{2B}) \mbox{   by Lemma \ref{LCox}} 
 \end{eqnarray*}
 establishing Proposition \ref{P1234-hack}.

\subsection{Size-biasing the percolation process and Yule process}
\label{sec-sizebias}
This section builds up to proving a result, Proposition \ref{prop:size-bias-growth}, about the number of vertices seen by the percolation process from vertex $1$ when we condition on existence of a short path of specified length from vertex $1$.

\subsubsection{Some terminology}
Let us quickly revisit the structures (section \ref{sec:prelim}) within $\GG_n$ associated with percolation from vertex $1$ and introduce more precise terminology.  
The {\em percolation tree} itself is the spanning tree consisting of all edges in the shortest paths
$\bpi(1,j), \ 2 \leq j \leq n$.  
The {\em percolation tree process} tells us at time (time = distance) $t$ the subtree on vertices within distance $t$ from vertex $1$.
And the {\em percolation counting process} $N_n(t)$ at (\ref{def-Nnt}) tells us at time $t$ the number of vertices within distance $t$ from vertex $1$.
We can use the same terminology for the Yule process of section \ref{sec-Yule}; the process 
$(N_\infty(t), \ 0 \leq t < \infty)$ is the {\em Yule counting process}.  
The underlying continuous-time branching process, run until time $t$ and then regarded as a random tree with edge-legths, is the {\em Yule tree process} at time $t$.  
This process run to time $\infty$ is the {\em Yule tree} or {\em PWIT}, a random infinite tree with edge-lengths.

\subsubsection{Heuristics for size-biasing}
Associated with the Yule counting process is the limit 
(Lemma \ref{LYule}(b)) 
random variable 
$W:= \lim_t e^{-t}N_\infty(t)$ 
with  probability density $e^{-w}$ on $ 0 < w < \infty$.  
What can we say about the Yule tree conditioned on it having a vertex at some specified large distance 
$t_0$?  The probability of this event given $W$ is approximately proportional to $W$, so the posterior density  of $W$ given this event becomes approximately $we^{-w}$.   
This is a basic instance of size-biasing.
But instead of relying  on Bayes calculations for single random variables, we describe next the more elegant approach to size-biasing the whole process  based on a probabilistic construction (this type of construction is widely used in  modern branching process theory
\cite{LPP95c}).  
In this method the density $we^{-w}$ arises as the density of the sum of two independent Exponential($1$) random variables.

\subsubsection{The size-biased Yule process}
\label{sec-sbY}
We are working toward a result of the type 
\begin{quote}
the $n \to \infty$ limit of the size-biased percolation process is the size-biased Yule process
\end{quote}
and now we will define and derive simple properties of the limit process, without justifying the
``size-biased" name.

\begin{figure}
\begin{center}

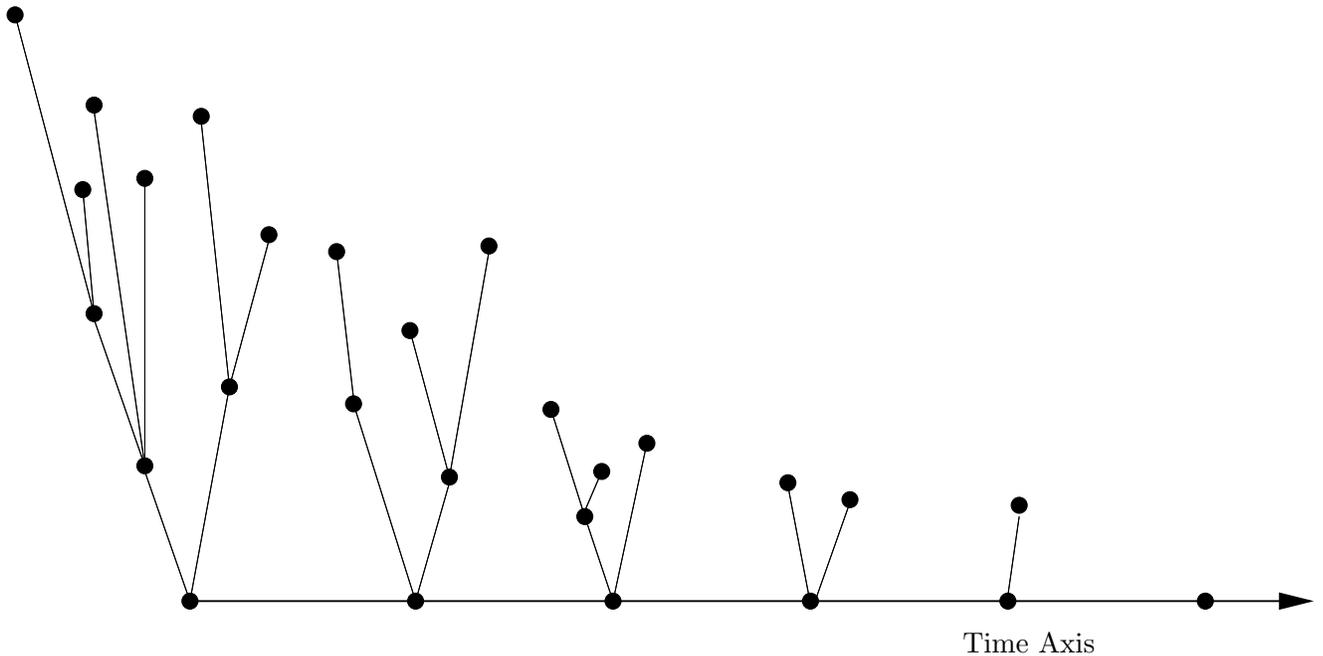

\caption{{\bf Size-biased Yule process}} \label{size-bias}
\end{center}
\end{figure}

On the half-line $\Rbold^+$, put a ``root" vertex at the origin and other vertices at the points $(P_i, i \geq 1)$ 
 of a rate 1 Poisson process.   
 Make each of these vertices the root of a Yule tree.  
 Regarding the resulting structure as a random infinite tree with edge-lengths, call it the {\em size-biased Yule tree}, with root at the origin.  
 Given a distance $t$, the {\em size-biased Yule process} at $t$ is the subtree on vertices at distance less than $t$ from the root, illustrated in Figure  \ref{size-bias}.  
 The associated counting process, giving the number  of vertices at distance less than $t$ from the root, 
 is
 \begin{equation} \widetilde{N}(t) = N_0(t) + \sum_{i:P_i \leq t} N_i(t-P_i) \label{Ntilde} 
\end{equation}
 where $(N_0,N_1,N_2,\ldots)$ are the Yule counting processes associated with the constituent Yule trees.
 
 Call the original half-line the \emph{distinguished path to infinity.}  
 We will also use, for technical reasons, the variation where the distinguished path is cut at some at some large finite distance $s$ from the origin, so its counting process is 
\[  \widetilde{N}^s(t)  =  N_0(t)+ \sum_{i: P_i \leq \min(s,t) } N_i(t-P_i) . \]

\vspace{0.2in}

We collect below some simple facts about these two processes. Our main aim is to understand the limiting behavior of $\widetilde{N}(t)$ for large t, and similarly, the  behavior of $\widetilde{N}^s(t)$ for large $s$ and $t$.
\begin{Lemma}
\label{sby-growth}
(a) There exists a random variable  $\widetilde{W}$ with probability  density  $w e^{-w}, \ w>0$ such that
\begin{equation}
\lim_{t \to \infty} e^{-t} \widetilde{N}(t)  =  \widetilde{W}; \quad 
\lim_{s\rightarrow \infty} \lim_{t \rightarrow \infty} e^{-t} \widetilde{N}^s(t) = \widetilde{W}
\label{tilde-as}
\end{equation} 
where the convergence holds a.s. and in $L^2$.\\
(b) For any $c, s, t > 0$
\[
\Ex \widetilde{N}^s(t) 
\left(c + \widetilde{N}^s(t) + \#\{i: t \leq P_i \leq s\} \right)
\leq 
6e^{2t} + 2e^t(c+s) . \]
\end{Lemma}
\proof
The construction (\ref{Ntilde}) implies
 that the size biased process $\widetilde{N}(\cdot)$ can be represented in terms of the sum of two independent Yule processes
\begin{equation} \widetilde{N}(t) = N_0(t) + (N^\prime(t)-1)\label{NNN}
\end{equation}
 because the contribution from the distinguished path to infinity behaves as another Yule process rooted at the origin, with the distinguished path representing the reproduction times of the initial ancestor.  We subtract 1 to avoid double counting the root. 
Lemma \ref{LYule}(b) says we have independent Exponential($1$) limits (a.s. and in $L^2$)
\[ W_0 := \lim_t e^{-t} N_0(t);  \quad 
 W^\prime = \lim_t e^{-t} N^\prime(t) \] 
and (\ref{tilde-as}) follows easily.
For (b), because 
$\widetilde{N}^s(t)$ is independent of $\#\{i: t \leq P_i \leq s\}$,
the quantity under consideration equals
\[ (c + (s-t)^+) \Ex \widetilde{N}^s(t) + \Ex (\widetilde{N}^s(t))^2 . \]
Use the inequality 
$\widetilde{N}^s(t) \leq \widetilde{N}(t)$ 
and the inequalities (from (\ref{NNN}) and Lemma \ref{LYule}(a))
\[ \Ex \widetilde{N}(t) \leq 2e^t; \quad 
\var \widetilde{N}(t) \leq 2e^{2t} \]
to complete the proof of (b).
\qed

\subsubsection{The percolation counting process conditioned on existence of a path}
\label{growth-rate-one-path}
Proposition \ref{prop:size-bias-growth} will 
formalize the idea
\begin{quote}
Conditional on existence of a path from vertex $1$ of specified length, 
the percolation process is approximately the size-biased Yule process.
\end{quote}
In the following lemma,  ``number of vertices" excludes vertex $1$, and
we are conditioning on length-to-boundary being $s$.
\begin{Lemma}
\label{cond-exact}
Fix $\alpha \in L(\bt)$. Condition on the existence of a path of length $s + b(\alpha)$ from vertex 1 to vertex $\alpha$. Let $\widetilde{P}^{s}_n$ denote the number of vertices on this path, and $U_{(1)} < U_{(2)} < \ldots < U_{(\widetilde{P}_n^{s})}$ denote the distances of these vertices from vertex $1$. Then\\
(a) The exact distribution of $\widetilde{P}_n^{s}$ is 
\begin{equation}
\label{eq:exact-count}
\prob(\widetilde{P}_n^{s} = k ) = C(s,n)\frac{(n-1-\#\bt)_k}{n^k} \frac{s^k}{k!} 
, \quad 0 \leq k \leq n -1-\#\bt  
\end{equation}
where $C(s,n) = \left( \sum_{j=0}^{n-1-\#\bt} \frac{(n-1-\#\bt)_j}{n^j} \frac{s^j}{j!}\right)^{-1}$ is the normalizing constant.\\
(b) Conditional on $\widetilde{P}_n^s$ the % $U_i$'s are constructed as follows:\\
$(U_{(k)})$ are distributed as the order statistics of 
$\widetilde{P}_n^s$ independent Uniform$(0,s)$ random variables.
\\(c) Suppose $s_n \to \infty, \ s_n = o(\sqrt{n})$.
Then the variation distance between 
the distribution of $\widetilde{P}_n^{s_n}$ and the Poisson($s_n$) distribution 
tends to $0$ as $n \to \infty$.
\end{Lemma}
\proof
Formula (\ref{intens-1}) says that the intensity of the event 
\begin{quote}
[[there exists a path of length $s+b(\alpha)$ from $1$ to $\alpha$ whose vertices are at distances 
$0<u_1< \ldots < u_j < s$ ]]
\end{quote}
is of the form 
$c(n,s,\alpha) \ 
(n-1-\#\bt)_j n^{-j} 
$.
Now (a) follows from the integral identity (\ref{kfold}), 
and (b) follows from the uniformity of the density in $(u_1,\ldots,u_j)$.
And (c) follows from (a) because the ratio 
$\frac{(n-1-\#\bt)_j}{n^j}$
tends to $1$ for $j = o(\sqrt{n})$. 
\qed

For the main result of this section, we study a certain 
{\em pruned} percolation process  which we now define carefully.  
Recall that the percolation counting process $N_n(t)$ counts the number of vertices $j$ 
such that there exists a path $\pi$  from $1$ to $j$ of length $\leq t$.  
For the {\em pruned} percolation counting process we impose two extra restrictions on $\pi$. 
First, $\pi$ must not use any vertex in the neighborhood 
$\NN_\tau(n-1,n) = \bt_{A_n}$.
Next, we will be conditioning on existence of a path, say $(1,\eta_1,\eta_2,\ldots)$, 
of specified length from vertex $1$.  
Say $\pi$ contains a {\em short-cut} if  the path $\pi$ meets $\eta_j$  for some $j \geq 1$.
The second restriction on $\pi$ is that $\pi$ must not contain any short-cuts.

\begin{Proposition}
 \label{prop:size-bias-growth}
 Consider a sequence $s_n$ satisfying $\omega_n\leq s_n\leq \log{n}+B$, and a vertex $\alpha \in L(\bt)$. Condition on the existence of a path from $1$ to $\alpha$ of length $s_n+ b(\alpha)$. Let ${\widetilde N}^{s_n}_n(t)$ be the pruned percolation counting process defined above. Then for each $n$ there exists a  random variable ${\widetilde W}_n$  having density $we^{-w}$ on $\Rbold^+$ such that 
  \[\sup_{\omega_n \leq t\leq \log{n} -\omega_n}  \left| e^{-t} {\widetilde N}_n^{s_n}(t) - W_n\right| \longrightarrow 0 \mbox{ in probability} \]
as $n\to \infty$.
 \end{Proposition}
 We give the proof in some detail; later (Proposition \ref{prop:growth-four-sources}) we need the variant for percolation from several sources, and we will omit details of that variant.
 
\subsubsection{Proof  of Proposition \ref{prop:size-bias-growth}}
We start with some finite error bounds for the two size-biased Yule processes ${\widetilde N}(\cdot)$ and ${\widetilde N}^s(\cdot)$ that were introduced in section \ref{sec-sbY}.

\begin{Lemma}
 \label{ld}
 Consider any sequences $\omega_n, s_n  \to \infty$. 
 \\(a) Fix  $\eps> 0$. Recall the limiting random variable ${\widetilde W}$ from Lemma \ref{sby-growth}. Then there exists a constant $C$ such that
 \[\prob\left(\sup_{t\geq \omega_n} \left| e^{-t} {\widetilde N}(t) -{\widetilde W}  \right| > \eps \right) \leq C e^{-\omega_n}\cdot \eps^{-2} .\]  
 (b) Fix  $\eps> 0$ and consider the processes ${\widetilde N}^{s_n}(\cdot)$. Then there exist random variables ${\widetilde W}_n$ having density $we^{-w}$ on $\Rbold^+$ such that 
 \[\prob\left( \sup_{t \geq \omega_n} \left| e^{-t} {\widetilde N}^{s_n}(t) - {\widetilde W}_n \right| > \eps  \right)\leq C \eps^{-2} \left(e^{-2\cdot s_n} + e^{-\omega_n} \right) .\]
 \end{Lemma}
 \proof 
 From (\ref{NNN}) we see that  $e^{-t} \cdot ({\widetilde N(t)}+1)$ is a martingale. 
 Part (a) follows from Lemma \ref{sby-growth} and  the  $L^2$ maximal inequality for martingales.
 
 To prove part(b), let $Y(t)$ be a Yule process independent of ${\widetilde N}^{s_n}(\cdot)$ and define the process 
  \begin{eqnarray*}
 Z_n(t) &\! \! \! \! \! \! \! \! \! \! \! \! \! = \ \ \ \ \ \ \ \ \ \ {\widetilde N}^{s_n}(t) & \mbox{ for } t <  s_n\\
 &\ \ \ \ \ = {\widetilde N}^{s_n}(t) + Y(t-s_n)-1 & \mbox{ for }  t\geq s_n .
 \end{eqnarray*}
 Note that the process $Z_n(\cdot)$ has the same distribution as the untruncated size biased Yule process ${\widetilde N}(\cdot)$. Thus  there exists a limiting random variable ${\widetilde W}_n$ with density $x\cdot e^{-x}$ such that  inequality (a) is satisfied with $Z_n(\cdot)$ in place of ${\widetilde N}(\cdot)$.  Now note that for any $t> s_n$
 
 \[e^{-t}Z_n(t) = e^{-t} {\widetilde N}^{s_n}(t) + e^{-s_n} \cdot \left[e^{-(t-s_n)}Y(t-s_n)\right]  - e^{-t} .\]  
 For our desired asymptotics we can ignore the final $- e^{-t}$ term, and write 
\[\prob\left( \sup_{t\geq \omega_n} \left|\frac{{\widetilde N}^{s_n}(t)}{e^{t}} - {\widetilde W}_n \right| > \eps  \right)
\leq \prob\left(\sup_{t\geq \omega_n} \left|\frac{Z_n(t)}{e^{t}} - W_n\right|> \frac{\eps}{2}\right) + \prob\left(e^{-s_n}\cdot \sup_{t\geq 0}\frac{Y(t)}{e^{t}} > \frac{\eps}{2} \right)
 .  \]
Applying the $L^2$ maximal inequality to the martingale $e^{-t}Y(t)$ gives
 $\prob\left(e^{-s_n}\cdot \sup_{t\geq 0}\frac{Y(t)}{e^{t}} > \frac{\eps}{2} \right) \leq C\eps^{-2} e^{-2\cdot s_n}$.
 Combine with part(a) of the Lemma applied to $Z_n(\cdot)$ to get the result. 
 \qed
 
 We now give a construction of the 
 pruned percolation counting process 
 ${\widetilde N}^{s_n}_n(t)$, 
 designed for comparison with a similar construction of the 
 size-biased Yule process.
 Recall we are conditioning on  existence of a ``distinguished" path from $1$ to $\alpha$ of length $s_n+ b(\alpha)$. 
 Write 
 $0 = U_{0} < U_{1} < U_{2} < \ldots < U_{P_n} < s_n$ for the distances from $1$ to the vertices within this path.  
 Define
 \[A_n(t) := \#\{i \geq 0: \ U_i \leq  t \} . \] 
 We can write 
 \begin{equation}
 {\widetilde N}_n^{s_n}(t) = A_n(t)+ G_n(t) 
 \label{NAG}
 \end{equation}
 where $G_n(t)$ is the number of vertices in the 
 pruned percolation counting process 
 which are not on the distinguished path.  
 By definition the process $G_n(\cdot)$ evolves as the counting process satisfying
 \begin{equation}
 \prob\left(\left. G_n(t+dt) - G_n(t) = 1 \right| \GG_n(t) \right) 
 = n^{-1} \ (A_n(t)+G_n(t)) \ (n-P_n - \# \bt - G_n(t)) \ dt 
 \label{Gn-def}
 \end{equation}
 because the number of vertices wetted at $t$ equals $A_n(t)+G_n(t)$ and the number of available vertices to be wetted  equals $n-P_n - \# \bt - G_n(t)$, the terms $P_n$ and $\#\bt$ arising from the two restrictions in the definition of {\em pruned}.  
 The filtration used here has $\GG_n(0)$ as the $\sigma$-field generated by the $(U_i)$ and then 
 $\GG_n(t) = \sigma(\GG_n(0), G_n(s), 0 \leq s \leq t)$.
 
 To relate this construction to the size-biased Yule process, 
 it does no harm to assume 
 (by variation distance convergence, Lemma \ref{cond-exact}(c)) 
 that $P_n$ has exactly Poisson($s_n$) distribution, so that 
 $(U_1,\ldots,U_{P_n})$ are the points of a rate-$1$ Poisson point process on $(0,s_n)$.  
 Now we could construction the size-biased Yule process, cut at $s_n$, via 
 \[  {\widetilde N}^{s_n}(t) = A_n(t)+ C_n(t) \] 
 where $C_n(\cdot)$ evolves as the counting process satisfying 
 \begin{equation}
   \prob\left(\left. C_n(t+dt) - C_n(t) = 1 \right| \GG(t) \right) 
 = (A_n(t)+C_n(t)) \ dt  
 \label{Cn-def}
 \end{equation}
 for appropriate filtration $(\GG(t))$.  
 But it is more useful to couple the two processes by first defining ${\widetilde N}_n^{s_n}(\cdot)$ 
 via (\ref{NAG}) and then defining ${\widetilde N}^{s_n}(\cdot)$ via
  \begin{equation}
 {\widetilde N}^{s_n}(t) = A_n(t)+ G_n(t) + B_n(t)
 \label{NAGB}
 \end{equation}
  where $B_n(\cdot)$ evolves as the counting process with $B_n(0) = 0$ and 
  \begin{equation}
    \prob\left(\left. B_n(t+dt) - B_n(t) = 1 \right| \GG_n(t) \right) 
 = (b_n(t) + B_n(t))\ dt  
 \label{Bn-def}
 \end{equation}
 where (subtracting (\ref{Gn-def}) from (\ref{Cn-def}))
  $b_n(t) + B_n(t) = (A_n(t) + G_n(t) + B_n(t)) -  n^{-1} \ (A_n(t)+G_n(t)) \ (n-P_n - \# \bt - G_n(t))$, 
 which works out as
 \[
 0 \leq b_n(t) = (A_n(t)+G_n(t)) \ \sfrac{P_n + \# \bt + G_n(t)}{n}  .
 \] 
 In particular, $ {\widetilde N}^{s_n}_n(t) \leq  {\widetilde N}^{s_n}(t)$, and $B_n(\cdot)$ is the number of extra vertices in the size-biased Yule process but not in the pruned percolation process.
 In view of Lemma \ref{ld}(b),  to prove Proposition \ref{prop:size-bias-growth} it is sufficient to prove
  \begin{equation}
  \sup_{\omega_n \leq t\leq \log{n} -\omega_n}   e^{-t} B_n(t)  \longrightarrow 0 \mbox{ in probability.} 
  \label{omega_n}
  \end{equation}
 Note that we can write 
 \begin{eqnarray*}
  b_n(t) &=& n^{-1} \  {\widetilde N}^{s_n}_n(t) \ ( {\widetilde N}^{s_n}_n(t) + \#\bt + \#\{i: t \leq U_i \leq s_n\}) \\
  &\leq&  n^{-1} \  {\widetilde N}^{s_n}(t) \ ( {\widetilde N}^{s_n}(t) + \#\bt + \#\{i: t \leq U_i \leq s_n\}) 
 \end{eqnarray*} 
 and then Lemma \ref{sby-growth}(c) implies 
   \begin{equation}
    Eb_n(t) \leq n^{-1}(6e^{2t} + 2e^t(\#\bt +s_n) ) 
 := a_n(t), \mbox{ say.} 
   \label{bnan}
  \end{equation}
 Consider the event
 \[ \Omega_n := \{B_n(\sfrac{1}{3} \log n) = 0\} \]
 that the two processes coincide up to time $\sfrac{1}{3} \log n$.  
 Using (\ref{Bn-def}) 
 \[ 1 - \prob (\Omega_n) \leq \int_0^{\sfrac{1}{3} \log n} \Ex b_n(t) dt 
 \leq \int_0^{\sfrac{1}{3} \log n} a_n(t) dt \to 0 \]
 and so $\prob (\Omega_n) \to 1$.  
 Next observe that 
 $e^{-t}B_n(t)$ is a submartingale, because 
 \[ d(e^{-t}B_n(t)) = e^{-t} dB_n(t) - e^{-t} B_n(t) dt \]
 and so 
 \[ \Ex( d(e^{-t}B_n(t)) | \GG_n(t)) = e^{-t} b_n(t) dt \geq 0 . \] 
 Appealing to the $L^1$ maximal inequality for submartingales, to prove (\ref{omega_n}) it is now enough to prove
  \begin{equation}
 \Ex  e^{-t_n} B_n(t_n) \ind(\Omega_n)   \longrightarrow 0 \mbox{ for } t_n = \log n - \omega_n .
  \label{omega_z}
  \end{equation}
  Write 
  $f_n(t) =  \Ex  B_n(t) \ind(\Omega_n)$ for $t \geq \sfrac{1}{3} \log n$, so that 
  $f_n(\sfrac{1}{3} \log n) = 0$.
  Using (\ref{Bn-def}), 
  $ f^\prime_n(t) \leq f_n(t) + a_n(t)$ and so 
  \[ (e^{-t}f_n(t))^\prime \leq e^{-t} a_n(t) . \]  
  Using (\ref{bnan}), for $t \geq \sfrac{1}{3} \log n$ we have 
  $a_n(t) \leq C e^{2t}/n$ for some constant $C$.
  So 
  $(e^{-t}f_n(t))^\prime \leq C e^t/n$ and then (\ref{omega_z}) holds because
  \[ e^{-t_n}f_n(t_n) = \int_{\sfrac{1}{3} \log n}^{t_n} C e^t/n \ dt 
  \leq Cn^{-1} \exp(t_n) \to 0 . \]

\subsection{Proof of Proposition \ref{S1234}}
\label{sec-P22}
Proposition \ref{prop:size-bias-growth}
studied the pruned percolation counting process starting from vertex $1$.
We now want to consider  four such processes running 
concurrently, starting from vertices $1,2,3,4$ (``sources").
In this setting, if one of the flows reaches a vertex $j$ 
which was previously reached by a different flow,
we say a {\em collision} occurs, and vertex $j$ is only counted in the counting process ($\widetilde{N}_n^{(i)}(t)$ below)
for the source $i$ whose flow first reaches $j$.

Recall the setting of Proposition \ref{S1234}: 
we say ``conditional on $s_1,s_2,t_1,t_2$" to mean conditional on the event described at the beginning of section \ref{cond-distr-other-short}.  Note that  we condition only on the \emph{lengths} and not 
on the internal structure of the four distinguished path segments.  
The values of $s_1,s_2,t_1,t_2$ (which depend on $n$) are assumed to satisfy
\begin{equation}
\omega_n \leq s^1 := s_1 + s_2 + \sigma \leq \log n + B; \quad 
\omega_n \leq t^1 := t_1 + t_2 + \sigma \leq \log n + B. 
\end{equation}
For $i = 1,2,3,4$ write 
$\widetilde{N}_n^{(i)}(t)$ for the number of vertices reached by the flow started at source $i$ before time t, in the {\em concurrent flow process}.  
This process differs from the $4$ separate processes in two ways.
First, the elimination of collisions, as described above. 
Second, we extend the notion of (forbidden) short-cuts to say that a path of the percolation process may not meet any of the $4$ distinguished paths.
Of course for each $i$, the number $\widetilde{N}_n^{(i)}(t)$ is bounded by the corresponding number 
in the percolation flow from $i$ when the other flows are not present, 
which was the context of Proposition \ref{prop:size-bias-growth}.

\begin{Proposition}
\label{prop:growth-four-sources}
Conditional on  $s_1,s_2,t_1,t_2$, 
there exist random variables 
$\widetilde{W}^{(i)}_n$ such that
\[ (\widetilde{W}^{(1)}_n, \widetilde{W}^{(2)}_n, \widetilde{W}^{(3)}_n, \widetilde{W}^{(4)}_n) \cd  
(\widetilde{W}^{(1)}, \widetilde{W}^{(2)}, \widetilde{W}^{(3)}, \widetilde{W}^{(4)}) \]
where the limit r.v.'s are independent with density $we^{-w}$; 
and such that, for any 
$\omega_n \leq t_n \leq \frac{1}{2}(\log n + B)$,
\begin{equation}
 e^{-t_n}  \widetilde{N}^{(i)}_n(t_n)  - \widetilde{W}^{(i)}_n
\to 0 \mbox{ in probability} 
\label{Nst2}
\end{equation}
for each $1 \leq i \leq 4$.
\end{Proposition} 
\proof
The proof involves only minor modifications of the proof of Proposition \ref{prop:size-bias-growth} -- the essential issue is to show that the two changes (collisions; short cuts) in going from separate to concurrent processes has negligible effect.   We omit details.
\qed

Recall the definition of $S_{12},S_{34} $ and the definition of the Cox point processes from section \ref{cond-distr-other-short}. Let $S^{*}_{12},$ (resp. $S^{*}_{34} $)  be the times of the first collision (within the  concurrent flow process) between the flow processes starting from 1 and 2 (resp. from 3 and 4).  Note that if a collision occurs between the flow processes started at 1 and 2 at time $t$, then there is a path of length $2t$ from 1 to 2, and 
(because we do not allow flow through the neighborhood $\NN_\tau(n-1,n)$) 
this path does not use the distinguished edge of the neighborhood.  
So $S_{12} \leq 2 S^{*}_{12} $ and $S_{34} \leq 2 S^{*}_{34} $.  
This is an {\em in}equality because there might be shorter paths that were ``pruned away" in the processes we have studied.  So  
to prove Proposition \ref{S1234} it is enough to prove the following
\begin{Proposition}
Conditional on  $s_1,s_2,t_1,t_2$, 
\begin{equation}
\label{converge-law}
(2S^{*}_{12} -\log{n},   2 S^{*}_{34} -\log{n})  \cd( \xi_1^{12},  \xi_1^{34} )
 \end{equation}
as $n \rightarrow \infty $.
\end{Proposition}
\proof 
Condition on the concurrent flow process until time $t$, and suppose the flows from source $1$ and source $2$ have not collided before time $t$.  
Then the instantaneous conditional probability-per-unit-time of a collision (``hazard rate") equals 
\begin{eqnarray*}
\lambda_n(t) &=& \frac{2 \widetilde{N}_{n}^{(1)}(t) \ \widetilde{N}_{n}^{(2)}(t)}{n} 
\end{eqnarray*}
because the unseen length of each possible edge has Exponential($1/n$) distribution.  
We are interested in the recentered process $2S_{12}^{*} -  \log{n} $.  This process has hazard rate 
\begin{equation}
\widetilde{\lambda}_n(s)= \sfrac{1}{2} \lambda_n(\sfrac{1}{2}s + \sfrac{1}{2} \log n) =
e^{s} \  \frac{\widetilde{N}^{(1)}_{n}(\frac{1}{2}\log{n} + \frac{1}{2}s)}{e^{\frac{1}{2}\log{n} + \frac{1}{2}s}}.
\frac{\widetilde{N}^{(2)}_{n}(\frac{1}{2}\log{n} + \frac{1}{2} s)}{e^{\frac{1}{2}\log{n} +  \frac{1}{2}s}}
\end{equation}
Now use Proposition \ref{prop:growth-four-sources} to conclude
\[ e^{-s} \widetilde{\lambda}_n(s) \cp  \widetilde{W}^{(1)} \widetilde{W}^{(2)}
\mbox{ uniformly on }  2 \omega_n - \log n \leq s \leq B .
 \] 
 This easily implies 
$ 2S^{*}_{12} -  \log{n} \cd \xi_1^{12}$, 
because $\xi_1^{12}$ is defined to have hazard rate 
$e^s \widetilde{W}^{(1)} \widetilde{W}^{(2)}$ on $- \infty < s < \infty$.  
The joint convergence (\ref{converge-law}) follows by the same argument, the independence of the limits $(\widetilde{W}^{(i)} )$ in Proposition \ref{prop:growth-four-sources} implying independence of the limits 
$( \xi_1^{12},  \xi_1^{34} )$ here.
\qed.

\section{Further discussion}
\label{sec-FD}
\subsection{Analysis of the limit function $G(z)$}
\label{sec-Gz}

Recall
$$ G(z)=\int_{0}^{\infty}P(W_1W_2e^{-u}>z)\ du $$
where $W_1$ and $W_2$ are independent
Exponential($1$).
The Mellin transform of % the given function is given by
$G$ is
%\begin{displaymath}
%\parbox{.75\textwidth}{
\begin{eqnarray*}
\Phi(y)&=& \int_0^{\infty}z^{y-1}G(z)dz \\
&=&  \int_{0}^{\infty}\int_{0}^{\infty}z^{y-1}P\left(W_1>\frac{z e^u}{W_2}\right) \ dudz\\
&=&
\int_{0}^{\infty}\int_{0}^{\infty}\int_{0}^{\infty}\exp(-w)z^{y-1}\exp\left(-\frac{z{e^u}}{w}\right) \ dzdudw\\
%&=&\int_{0}^{\infty}\int_{0}^{\infty}e^{-(w)}\int_{0}^{\infty}z^{y-1}e^{-(\frac{e^u}{w}z)}
%dzdudw\\
&=&\int_{0}^{\infty}\int_{0}^{\infty}e^{-w}\Gamma(y)\left({\frac{e^u}{w}}\right)^{-y} \ dwdu\\
&=&
\Gamma(y)\int_{0}^{\infty}e^{-uy}\int_{0}^{\infty}e^{-w}w^y
dwdu\\
&=&\Gamma(y)\Gamma(y+1)\int_{0}^{\infty}e^{-uy}du\\
&=& \frac{\Gamma(y)\Gamma(y+1)}{y}\\
&=& {(\Gamma(y))}^2 . \\
\end{eqnarray*}
%}
%\end{displaymath}
Checking a table of Mellin transforms (\cite{MR0352890} II.5.34) we see 
\[ G(z) = 2 K_0(2 z^{1/2}) \]
where $K_0$ is the modified Bessel function of the second kind. 
The standard asymptotics of $K_0$ (\cite{MR1688958} 4.12.6) say
\[ K_{0}(x)  \sim \sqrt \frac{\pi}{2x}\exp{(-x)} \quad \mbox{ as } x \to \infty \]
and so
 \[ G(z) \sim \pi^{1/2} z^{-1/4} \exp{(-2\sqrt{z} )} \quad  \mbox{ as } z \to \infty . \]

\begin{figure}

\begin{center}

\centerline{\epsfig{figure=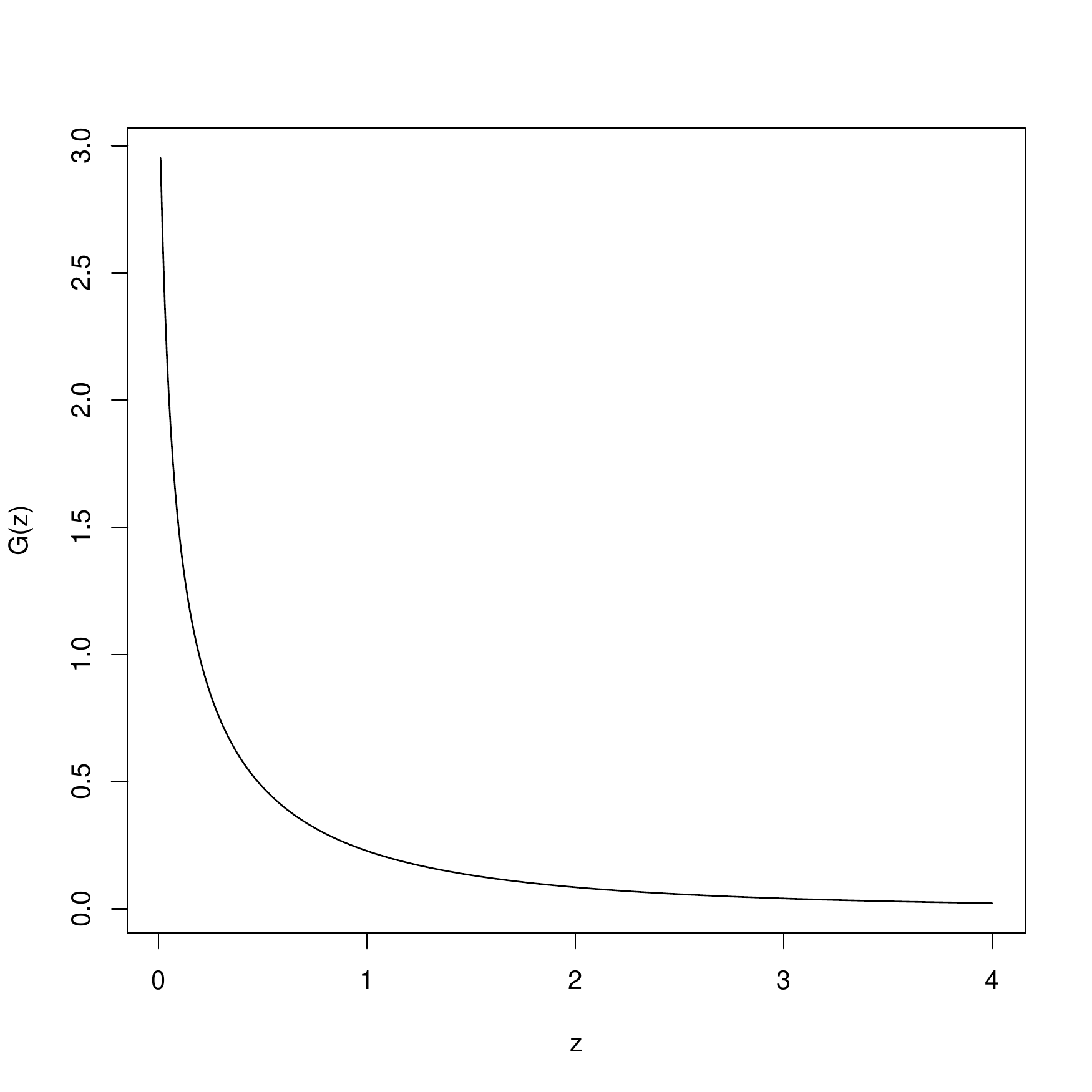,height=200pt}}

\end{center}

\caption{The function $G(z)$, drawn using Mathematica's numerical
 integration toolbox.}
\end{figure}

\subsection{Methodology of relating local and global structure}
\label{sec-meth}
As illustrated by the heuristic argument in section \ref{sec-heur},
the conceptual point of Theorem \ref{T1} is that a quantity
depending on the ``global" structure of the network can be studied
statistically via a ``local" (i.e. large fixed distance)
calculation.
This reduction to local structure is, to our understanding,
the central point in the powerful non-rigorous
{\em cavity method} of statistical physics \cite{MParisi03}.
In our attempted mathematical reformulations of the cavity method as applied
to combinatorial optimization problems such as TSP \cite{me101,me103,me109}
in this random network model $\GG_n$,
we made explicit use of the $n \to \infty$ limit structure (the PWIT of section \ref{sec-Yule})
of this model as viewed from a random vertex.
In these harder problems one needs rather abstract, often as yet not rigorously
justified, arguments to connect local and global structure.
The problem in this paper seems conceptually easier
in that we can use concrete calculations instead.

\subsection{Flows through vertices}
\label{sec-vertexflows}
In the setting of Theorem \ref{T1} 
one could alternatively consider flows through 
{\em vertices} instead of edges.
Let us state this alternative result and indicate the derivation
of the limit distribution without giving details of proof.

Write $F_n^*(v)$ for the flow through vertex $v \in [n]$.
Let $(W_i, i \geq 1)$ be independent Exponential($1$) r.v.'s and let
$0 < \xi_1 < \xi_2 < \ldots $
be the points of a Poisson (rate $1$) process on $(0,\infty)$.
Define
\[ \Xi : = \sum_i \sum_{j \neq i} W_iW_j \exp(-\xi_i - \xi_j) . \]
\begin{Corollary}
In the setting of Theorem \ref{T1},
as $n \to \infty$ for fixed $z>0$,
\[ \sfrac{1}{n} \# \{v: F_n^*(v) \leq z \log n\}
\to_{L_1} \Pr(\Xi \leq z) .\]
\end{Corollary}
The formula is most neatly derived using the $n \to \infty$
limit PWIT structure of $\GG_n$ 
\cite{me101}.
Relative to a typical vertex $v$ of the PWIT, the edge-lengths $(\xi_i)$ to adjacent vertices 
$(v_i)$ are distributed as the points of a Poisson (rate $1$) process on $(0,\infty)$.   
For each $v_i$ let $N_i(t)$ be the number of vertices within distance $t$ from $v_i$ using  paths not via $v$.  Then $e^{-t}N_i(t) \to W_i$ for i.i.d. Exponential($1$) r.v.'s $(W_i)$.  
The relative volume of flow through the two edges $v_i \to v \to v_j$ will then be 
$W_iW_j \exp(-\xi_i - \xi_j)$ by the argument for Theorem \ref{T1}.

\subsection{Different models of random networks}
\label{sec-different}
The heuristic argument of section \ref{sec-heur}
can be carried over to a variety of random networks models.
For example, fix a degree distribution
$P(\Delta = i), \quad i \geq 1 $
with finite $2+ \eps$ moment.
There are several ways
(e.g. the ``configuration model")
to formalize the idea
of a $n$-vertex graph which is random subject to the constraint
that the $n \to \infty$ asymptotic degree distribution is $\Delta$.
Such models have local weak limits which are simple 
branching processes;
looking outwards from
a typical edge $e$, each end-vertex is the founder
of a Galton-Watson branching process with offspring distribution
\[ P(\Delta^* = i) = (i+1) P(\Delta = i+1)/E\Delta . \]  
Thus we expect the average vertex-vertex distance $D_n$ in such a random graph to behave as
\[ \Ex \bar{D}_n \sim \frac{\log n}{\log \Ex \Delta^*} . \]
See \cite{math.PR/0605414} for proofs for several models.  
Now make a random network by assigning independent random lengths
$\eta_e$ to edges $e$ 
(note that here we do not scale
edge-lengths with $n$).
Then the Galton-Watson process above becomes a general Markov
branching process in which individuals have $\Delta^*$ offspring at
independent ages $(\eta_i)$; the population size process
$N_\infty(t)$ has some Malthusian growth constant $\theta$
and some a.s. limit
$\exp(- \theta t) N_\infty(t) \to Z$.
The heuristic argument from section \ref{sec-heur} now suggests
that the limit joint distribution of edge-lengths and relative edge-flows will be
\[ (\eta, Z_1 Z_2 \exp(- \theta \eta)) \]
where $(\eta, Z_1, Z_2)$ are independent.  
But giving a rigorous proof for the models in \cite{math.PR/0605414} may be technically challenging.

Finally, one might consider
models on the two-dimensional lattice with i.i.d. random edge-lengths.  
Here, studying {\em lengths} of shortest routes is tantamount to studying (unoriented) first passage percolation 
\cite{MR513421}.  
However, if $i$ is close to $i^\prime$ and $j$ is close to $j^\prime$ then we expect the routes 
from $i$ to $j$ and from $i^\prime$ to $j^\prime$ to coincide except near the endpoints.  
This suggests a quite different distribution of edge-flows, more specifically that $F(e)$ should have a power-law tail.

\subsection{Random demands}
\label{demands}
A small variation of our model is to assume that the total flow to be routed
from vertex $i$ to vertex $j$
is a random variable $D_{ij}/n$ instead of $1/n$; the flow is still routed along the same shortest path as in the uniform demand case.
So the flow across edge $e$ is
\[F_n(e) = \frac{1}{n} \sum_{i\in [n]} \sum_{j \in [n], j\neq i} D_{ij} 1\{e \in \bpi(i,j)\} .\]
Because the flow across an edge $e$ is made up from many different 
source-destination pairs, it is straightforward to 
add a ``law of large numbers" step to the proof of Theorem \ref{T1} 
and obtain the following corollary.
\begin{Corollary}
(a)  Suppose $D_{ij} \geq 0$ are independent with common mean $0 < \mu < \infty$ and with uniformly bounded second moments. Then
 \[\frac{1}{n}\#\{e: F_n(e) > z \mu \log{n}\} \rightarrow_{L^1} G(z) , \quad z > 0.\]
(b)  Suppose instead the {\bf gravitational model} $D_{ij} = D_i D_j$ where  $D_i \geq 0$ are  independent random variables with common mean $0 <\mu < \infty$ and with uniformly bounded second moments. Then 
 \[\frac{1}{n}\#\{e: F_n(e) > z \mu^2 \log{n}\} \rightarrow_{L^1} G(z), \quad z > 0 . \]
\end{Corollary}

\subsection{Joint distributions for shortest paths}
\label{sec-k-paths}
As described in section \ref{sec:prelim}, various aspects of shortest paths in the model $\GG_n$ have been studied.  The following ideas will be developed elsewhere. 
There is a known (implicitly, at least) limit distribution 
\[ D_n(1,2) - \log n \cd D(1,2) \]
for distance between a typical pair of vertices.  
Now fix $k \geq 3$. 
We expect a joint limit 
\begin{equation}
(D_n(1,2) - \log n,\ldots,D_n(1,k) - \log n) \cd (D(1,2), \ldots, D(1,k)) 
\end{equation}
and it turns out the limit distribution is 
\[ (D(1,2), \ldots, D(1,k)) \ed (\xi_1 + \eta_{12}, \ldots, \xi_1 + \eta_{1k}) \]
where $\xi_1$ has the double exponential distribution 
\[ \Pr (\xi \leq x) = \exp(-e^{-x}), \ - \infty < x < \infty \]
the $\eta_{1j}$ have logistic distribution 
\[ \Pr (\eta \leq x) = \sfrac{e^x}{1+e^x},  \ - \infty < x < \infty \] 
and (here and below) the r.v.'s in the limits are independent.  
Now we can go one step further: we expect a joint limit for the array
\[ (D_n(i,j) - \log n, 1 \leq i < j \leq k) \cd (D(i,j), 1 \leq i < j \leq k) \]
and the joint distribution of the limit is 
\[ (D(i,j), 1 \leq i < j \leq k) \ed (\xi_i + \xi_j - \xi_{ij},  1 \leq i < j \leq k)  \] 
where the limit r.v.'s all have the double exponential distribution.
This implies two representations for the original limit distribution:
\[ D(1,2) \ed \xi_1 + \eta_{12} \ed \xi_1 + \xi_2 - \xi_{12} . \]

%\bibliographystyle{plain}
%\bibliography{../../trees/alg,../../trees/networks,../../trees/me,../../trees/trees,../../trees/misc,../../trees/rwgbook}

\end{document}